\documentclass[12pt,leqno]{article}
\title{ $\pi$-adic approach of $p$-class group and unit group of $p$-cyclotomic fields}
\author{Roland Qu\^eme}
\usepackage{amsmath,amsthm,amstext,amsbsy}
\newtheorem{thm}{Theorem}[section]
\newtheorem{cor}[thm]{Corollary}

\newtheorem{lem}[thm]{Lemma}
\font\mathbb=msbm10
\newcommand{\N}{\mbox{\mathbb N}}
\newcommand{\R}{\mbox{\mathbb R}}

\newcommand{\Q}{\mbox{\mathbb Q}}
\newcommand{\Z}{\mbox{\mathbb Z}}
\newcommand{\s}{\mathbf s}
\newcommand{\modu}{\ \mbox{mod}\ }
\newcommand{\be}{\begin{equation}}
\newcommand{\ee}{\end{equation}}
\date{2005 fev 04}
\begin{document}
\maketitle
\tableofcontents
\clearpage
\abstract
Roland Qu\^eme

13 avenue du ch\^ateau d'eau

31490 Brax

France

2005 fev 04

tel : 0033561067020

mailto: roland.queme@wanadoo.fr

home page: http://roland.queme.free.fr/

				   
********************************************************************

Let $p>2$ be a prime. Let $\Q(\zeta)$ be the $p$-cyclotomic field. Let $\Q(\zeta)^+$ be its maximal totally real subfield.
Let $\pi$ be the prime ideal of $\Q(\zeta)$ lying over $p$.  
This articles aims to describe some $\pi$-adic congruences characterising the structure of the 
$p$-class group and of the $p$-unit group of the fields  $\Q(\zeta)$ and $\Q(\zeta)^+$. For the unit group, this paper supplements the papers of D\'enes of 1954 and 1956.
A complete summarizing of the results obtained in the paper   follows in the   {\it Introduction section} \ref{s301071} 
from p. \pageref{s301071} to \pageref{s502051} . This paper is at elementary level. 
\endabstract
%
\clearpage
\section{Introduction}\label{s301071}
%
Let $p>2$ be a prime. Let $\Q(\zeta)$ be the $p$-cyclotomic fields. Let $\Z[\zeta]$ be the ring of integers of $\Q(\zeta)$. Let 
$\pi=(1-\zeta)\Z[\zeta]$ be the prime ideal of $\Q(\zeta)$ lying over $p$.  
This monograph contains two parts: 
\begin{enumerate}
\item
a description of $\pi$-adic congruences strongly connected to $p$-class group of $\Q(\zeta)$ and its structure.
\item
a description of $\pi$-adic congruences  on $p$-unit group of $\Q(\zeta)$.
\end{enumerate}
%
\subsection{ Some $\pi$-adic congruences connected to $p$-class group of cyclotomic field $\Q(\zeta)$}
This topic is studied in section \ref{s09121} p. \pageref{s09121} of this paper.
Let us give at first some definitions:
\begin{enumerate}
\item
Let $p$ be an odd prime. Let $\zeta$ be a root of the equation $X^{p-1}+X^{p-2}+\dots+X+1=0$.
Let $\Q(\zeta)$ be the $p$-cyclotomic field and $\Z[\zeta]$ be the ring of integers of $\Q(\zeta)$.
\item
Let $\sigma:\Q(\zeta)\rightarrow\Q(\zeta)$ be a $\Q$-isomorphism of the field $\Q(\zeta)$ generating the cyclic Galois group $G=Gal(\Q(\zeta)/\Q)$. There exists $u\in\N$,  primitive root $\modu p$, such that 
$\sigma(\zeta)=\zeta^u$.
\item
Let $C_p, C_p^+,C_p^-$ be respectively the subgroups of exponent $p$ of the $p$-class group of $\Q(\zeta)$, of the $p$-class group of 
$\Q(\zeta+\zeta^{-1})$ and  the relative $p$-class group $C_p^-=C_p/C_p^+$.
Let $r_p,r_p^+,r_p^-$ be respectively the $p$-rank of $C_p, C_p^+, C_p^-$, seen as ${\bf F}_p[G]$ modules.
\item
It is possible to write $C_p$ in the form  $C_p=\oplus_{i=1}^{r_p} \Gamma_i$, where $\Gamma_i$ is a cyclic group of order $p$, subgroup globally invariant under the action of the Galois group $G$.
\item 
Let $\mathbf b_i,\quad i=1,\dots,r_p$, be  a not principal integral ideal of $\Q(\zeta)$ whose class belongs to  the group $\Gamma_i$. Observe at first that $\mathbf b_i^p$ is principal and that $\sigma(\mathbf b_i)\simeq \mathbf b_i^{\mu_i}$ where $\simeq$ is notation for class equivalence and $\mu_i\in{\bf F}_p^*$ with  ${\bf F}_p^*$  the set of $p-1$ no null elements of  the finite field of cardinal $p$.
Let the ideal $\mathbf b=\prod_{i=1}^{r_p}\mathbf b_i$, which generates $C_p$ under action of the group $G$. 
\item
A number $a\in \Q(\zeta)$ is said {\it singular} if there exists a not principal  ideal $\mathbf a$ of $\Q(\zeta)$ such that 
$a\Z[\zeta]=\mathbf a^p$.
A singular number $a\in \Q(\zeta)$ is said {\it primary} if there exists $\alpha\in \N, \quad \alpha\not\equiv 0\modu p$ such that 
$a\equiv\alpha^p\modu\pi^p$.
\item
Let $d\in\N,\quad p-1\equiv 0\modu d$. 
Let $G_d$ be the subgroup of order $\frac{p-1}{d}$ of the Galois group $G$.
Let us define the minimal polynomial $P_{r_d}(X)$ of degree $r_d$ in the indeterminate $X$, where $P_{r_d}(\sigma^d)\in{\bf F}_p[G_d]$  annihilates the ideal class of $\mathbf b$, written also  
$\mathbf b^{P_{r_d}(\sigma^d)}\simeq \Z[\zeta]$. The polynomial $P_{r_d}(\sigma^d)$ is of form
$P_{r_d}(\sigma^d)=\prod_{i=1}^{r_d}(\sigma^d-\mu_i^d),\quad \mu_i\in{\bf F}_p^*$.
When $d=1$, then $G_d=G$,  Galois group of $\Q(\zeta)/\Q)$,  and $r_d=r_1$. Note that if $d>1$ then $r_d\leq r_1$.
\end{enumerate}
We obtain the following results:
\begin{enumerate}
\item
For $d_1, d_2$ co-prime natural integers  with $d_1\times d_2=p-1$, the degrees in the indeterminate $X$ of minimal polynomials $P_{r_{d_1}}(X)$ and $P_{r_{d_2}}(X)$ verify:  if 
$r_{d_1}\geq 1$ and $r_{d_2}\geq 1$ then $r_{d_1}\times r_{d_2}\geq r_1$.
\item
Let us set $d=1$.
Let us note $r_1=r_1^++r_1^-$, where $r_1^+$ and $r_1^-$ are respectively the degrees of the minimal polynomials $P_{r_1^+}(\sigma)$ and $P_{r_1^-}(\sigma)$,   corresponding to annihilation of groups $C_p^+$ and $C_p^-$
with $C_p=C_p^+\oplus C_p^-$. The following result connects strongly the degree $r_1^-$ to Bernoulli Numbers:
the degree   $r_1^-$ is the index of irregularity of $\Q(\zeta)$ (the number of even Bernoulli Numbers $B_{p-1-2m}\equiv 0 \modu p$ for $1\leq m\leq \frac{p-3}{2})$. Moreover the degree $r_1^-$ verifies the inequality $r_p^--r_p^+\leq r_1^-\leq r_p^-$.
\item
Let the ideal $\pi=(\zeta-1)\Z[\zeta]$.
The following results are $\pi$-adic congruences strongly connected to structure of $p$-class group 
$C_p$ of $\Q(\zeta)$:
\begin{enumerate}
\item
There exists  singular algebraic integers 
$B_i\in\Z[\zeta]-\Z[\zeta]^*,\quad i=1,\dots,r_p$, verifying:
\begin{enumerate}
\item
$B_i\Z[\zeta]=\mathbf b_i^{p}$ with $\mathbf b_i$ defined above
\item
$\sigma(\mathbf b_i)\simeq\mathbf b_i^{\mu_i}$.
\item $\sigma(B_i)=B_i^{\mu_i}\times\alpha_i^p,\quad \alpha_i\in\Q(\zeta),\quad \mu_i\in{\bf F}_p^*$.
\item
$\sigma(B_i)\equiv B_i^{\mu_i} \modu\pi^p$.
\item 
For  the value $m_i\in\N$ verifying  $\mu_i=u^{m_i}\mod p,\quad  1\leq m_i\leq p-2$, then 
\begin{displaymath}
\pi^{m_i}\ |\  {B_i-1}.
\end{displaymath}
\end{enumerate}
\item
We can precise the previous result: with {it a certain reordering of indexing} of $B_i,\quad i=1,\dots,r_p$,  
\begin{enumerate}
\item
For $i=1,\dots,r_p^+$,  then the $B_i$ are  {\it primary}, so $\pi^p\  |\ B_i-1$.
\item
For $i=r_p^++1,\dots,r_p^-$, then the $B_i$ are  {\it not primary}. They verify the  congruence 
\begin{displaymath}
\pi^{m_i}\ \|\  B_i-1.
\end{displaymath}
\item
For $i=r_p^-+1,\dots,r_p$,  then the $B_i$ are {\it primary or not primary} 
(without being able to have a more precise result)
with 
\begin{displaymath} 
\pi^{m_i}\ |\ B_i-1.
\end{displaymath} 
\end{enumerate}
\item
Let $\mu_i=u^{2m_i+1}\modu p$ with $1\leq m_i\leq \frac{p-3}{2}$  corresponding  to an ideal $\mathbf b_i$ whose class belongs to $C_p^-$, relative $p$-class group of $\Q(\zeta)$. In that case define 
$C_i=\frac{B_i}{\overline{B}_i}$ with $B_i$ already defined, so with $C_i\in\Q(\zeta)$. If $2m_i+1>\frac{p-1}{2}$ then it is possible to prove the explicit very straightforward  formula for $C_i\modu \pi^{p-1}$: 
\begin{displaymath}
C_i\equiv 1-\frac{\gamma_{p-3}}{1-\mu_i} \times 
(\zeta+\mu_i^{-1}\zeta^u+\dots+\mu_i^{-(p-2)}\zeta^{u^{p-2}})\modu \pi^{p-1},\quad \gamma_{p-3}\in{\bf F}_p^*.
\end{displaymath}
\end{enumerate}
\end{enumerate}
%
\subsection
{Some $\pi$-adic congruences on $p$-unit group  the cyclotomic field}
This topic is studied in section \ref{s210163} p. \pageref{s210163}.
We apply in following results to unit group $\Z[\zeta+\zeta^{-1}]^*$ the method applied to $p$-class group in previous results: 
\begin{enumerate}
\item
There exists a fundamental system of units $\eta_i,\quad i=1,\dots,\frac{p-3}{2}$, of the group
$F=\{\Z[\zeta+\zeta^{-1}]^*/(\Z[\zeta+\zeta^{-1}]^*)^p\}/<-1>$ verifying the relations:
\begin{equation}\label{e201274}
\begin{split}
& \eta_i\in\Z[\zeta+\zeta^{-1}]^*,\quad i=1,\dots,\frac{p-3}{2},\\
& \sigma(\eta_i)=\eta_i^{\mu_i}\times\varepsilon_i^p,\\
&\varepsilon_i\in\Z[\zeta+\zeta^{-1}]^*,\\
&n_i\in\N,\ \mbox{with\ }\mu_i=u^{2n_i}\modu p,\quad 1\leq n_i\leq \frac{p-3}{2},\\
& \eta_i\equiv 1\modu \pi^{2n_i},\quad i=1,\dots,
\frac{p-3}{2},\\
& \sigma(\eta_i)\equiv \eta_i^{\mu_i}\modu\pi^{p+1},\quad i=1,\dots,\frac{p-3}{2}.
\end{split}
\end{equation}
\item
With a {\it certain reordering of indexing} of $i=1,\dots,\frac{p-3}{2}$,  
\begin{enumerate}
\item
For $i=1,\dots,r_p^+$ then  $\eta_i$ are {\it not primary units} and
\begin{displaymath}
\pi^{2n_i}\ \|\ \eta_i-1.
\end{displaymath}
\item
For $i=r_p^++1,\dots,r_p^-$, then $\eta_i$ are {\it primary units} and
\begin{displaymath}
\pi^{a_i(p-1)+2n_i}\ \|\ \eta_i-1,\quad a_i\in\N,\quad a_i>0.
\end{displaymath}
\item 
For $i=r_p^-+1,\dots,r_p$, then $\eta_i$ are {\it not primary or primary units}  and 
\begin{displaymath}
\pi^{a_i(p-1)+2n_i}\ \|\  \eta_i-1,\quad a_i\in\N,\quad a_i\geq 0.
\end{displaymath}
\item 
For $i=r_p+1,\dots,\frac{p-3}{2}$, then $\eta_i$ are {\it not primary units} and 
\begin{displaymath}
\pi^{2n_i}\ \| (\eta_i-1).
\end{displaymath}
\end{enumerate}
\item
If $2n_i>\frac{p-1}{2}$ then it is possible to prove the very straightforward  explicit formula
for $\eta_i$: 
\begin{displaymath}
\eta_i\equiv 1-\frac{\gamma_{p-3}}{1-\mu_i} \times 
(\zeta+\mu_i^{-1}\zeta^u+\dots+\mu_i^{-(p-2)}\zeta^{u^{p-2}})\modu \pi^{p-1},\quad \gamma_{p-3}\in{\bf F}_p^*.
\end{displaymath}
\end{enumerate}
\label{s502051}
%
%
%
%
%
%
\clearpage
\section{Cyclotomic Fields : some definitions} 
In this section, we fix notations used in all this paper.
\begin{itemize}
\item
For $a\in\R^+$, we note $[a]$ the integer part of $a$ or the integer immediately below $a$.
\item 
We denote $[a,b],\quad a,b\in \R$, the closed interval bounded by $a,b$.
\item 
Let us  denote $< a > $ the cyclic group generated by the element $a$.
\item 
Let $p\in\N$ be an odd prime.
\item
Let $\Q(\zeta_p)$, or more brievely $\Q(\zeta)$ when there is no ambiguity of the context, be the $p$-cyclotomic number field. 
\item
Let $\Z[\zeta]$ be the ring of integers of $\Q(\zeta)$.
\item
Let $\Z[\zeta]^*$ be the group of units of $\Z[\zeta]$.
\item 
Let $\Q(\zeta+\zeta^{-1})$ be the maximal real subfield of $\Q(\zeta)$, with
$[\Q(\zeta):\Q(\zeta+\zeta^{-1})]=2$. The ring of integers of $\Q(\zeta+\zeta^{-1})$ is $\Z[\zeta+\zeta^{-1}]$. Let $\Z[\zeta+\zeta^{-1}]^*$ be the group of units of 
$\Z[\zeta+\zeta^{-1}]$.
\item
Let ${\bf F}_p$ be the finite field with $p$ elements. 
Let ${\bf F}_p^*={\bf F}_p-\{0\}$.
\item 
Let us denote $\mathbf a$ the integral ideals of $\Z[\zeta]$. Let us note
$\mathbf a\simeq \mathbf b$ when the two ideals $\mathbf a$ and $\mathbf b$ are in the same class of the class group of $\Q(\zeta)$. The relation $\mathbf a\simeq\Z[\zeta]$ means that  the ideal $\mathbf a$ is principal.
\item
Let us  note $Cl(\mathbf a)$ the class of the ideal $\mathbf a$ in the class group of $\Q(\zeta)$.
Le us note $<Cl(\mathbf a)>$  the finite group generated by the class $Cl(\mathbf a)$. 
\item 
If $a\in\Z[\zeta]$, we note $a\Z[\zeta]$ the  principal integral ideal of $\Z[\zeta]$ generated by $a$.
\item 
We have $p\Z[\zeta]=\pi^{p-1}$ where $\pi$ is the principal prime ideal $(1-\zeta)\Z[\zeta]$.
Let us denote $\lambda= \zeta-1$, so $\pi=\lambda\Z[\zeta]$.
\item
Let $G=Gal(\Q(\zeta/\Q)$ be the Galois group of the field $\Q(\zeta)$.
Let $\sigma :\Q(\zeta)\rightarrow \Q(\zeta)$ be a $\Q(\zeta)$-isomorphism generating the cyclic group $G$. The $\Q$-isomorphism $\sigma$ can be defined by $\sigma(\zeta)=\zeta^u$ where $u$ is a primitive root $\modu p$. 
\item  
For this primitive root $u\modu p$ and $i\in\N$, 
let us denote $u_i\equiv u^i \modu p,\quad 1\leq u_i\leq p-1$.
For $i\in\Z,\quad i<0$, this is to be understood as $u_i u^{-i}\equiv 1 \modu p$.
This notation follows the convention adopted in Ribenboim \cite{rib}, last paragraph of page 118.
This notation is largely used in the sequel of this monograph.
\item
For $d\in\N,\quad p-1\equiv 0\modu d$, let $G_d$ be the cyclic subgroup of $G$ of order $\frac{p-1}{d}$ generated by $\sigma^d$, so with $G_1=G$. The group $G_d$ is the Galois group of the extension
$\Q(\zeta)/K_d$ where $K_d$ is a field with $\Q\subset K_d\subset \Q(\zeta)$ and $[K_d:\Q]=d$.
\item 
Let $C_p$ be the subgroup of exponent $p$ of the $p$-class group of the field $\Q(\zeta)$.
\item
Let $C_p^+$ be the subgroup of exponent $p$ of  the $p$-class group of the field $\Q(\zeta+\zeta^{-1})$. 
\item
Let $C_p^-$ be the relative class group defined by
$C_p^-=C_p/C_p^+$.
\item
Let  $h$ be the class number of $\Q(\zeta)$. The class number $h$ verifies the formula $h=h^-\times h^+$, where $h^+$ is the class number of the maximal real field $\Q(\zeta+\zeta^{-1})$, so called also second factor, and $h^-$  is the relative class number, so called first factor.
\item
Let us define respectively $e_p, e_p^-, e_p^+ $ by $h=p^{e_p}\times h_2,\quad h_2\not \equiv 0 \modu p$,
 by $h^-=p^{e_p^-}\times h_2^-,\quad h_2^-\not \equiv 0 \modu p$ and by 
$h^+=p^{e_p^+}\times h_2^+,\quad h_2^+\not \equiv 0 \modu p$
\item
Let $r_p, r_p^+, r_p^-,$ be respectively  the $p$-rank of the $p$-class group of $\Q(\zeta)$, of the $p$-class group of $\Q(\zeta+\zeta^{-1})$ and of the relative class group seen as ${\bf F}_p[G]$-modules, so with
$r_p\leq e_p,\quad r_p^+\leq e_p^+$ and $r_p^-\leq e_p^-$.
\item
The abelian group $C_p$ is a group of order $p^{r_p}$ with 
$C_p=\oplus_{i=1}^{r_p} C_i$ where $C_i$ are cyclic group of order $p$. 
\end{itemize}
%
\clearpage
\section{$\pi$-adic congruences on  $p$-subgroup $C_p$ of the class group of $\Q(\zeta)$}\label{s09121}
\begin{itemize}
\item
The two first subsections \ref{s20111} p.\pageref{s20111} and \ref{s20112} p.\pageref{s20112} give some definitions,  notations and general classical properties of the $p$-class group of the extension $\Q(\zeta)/\Q$. {\bf They can be omitted at  first}  and only looked at for fixing notations.
\item
In  subsection \ref{s108031} p. \pageref{s108031}, we get  several  results on  the structure of the $p$-class group $C_p$ of $\Q(\zeta)$ and on class number $h$ of $\Q(\zeta)$:
\begin{itemize}
\item
A formulation, with our notations,  of a Ribet's  result on irregularity index.
\item
Let $d,g\in\N$ coprime with $d\times g=p-1$. 
For groups generated by the action of  Galois groups $G$ and of subgroups 
$G_d, G_g$ of $G$  on ideals $\mathbf b$ of $\Q(\zeta)$,
an inequality between  degrees $r_1, r_d, r_g$ of minimal polynomials 
$P_{r_1}(\sigma)\in {\bf F}_p[G],\quad P_{r_d}(\sigma^d)\in{\bf F}_p[G_d],
\quad P_{r_g}(\sigma^g)\in {\bf F}_p[G_g]$  annihilating ideal class of $\mathbf b$.
\item
Some $\pi$-adic congruences  connected to structure of $p$-class group $C_p$ of $\Q(\zeta)$.
\end{itemize}
\end{itemize}
%
\subsection{Some definitions and notations}\label{s20111}\label{s111171}
In this subsection,  we fix or recall some notations used in all this section.
%
\begin{itemize}
\item
Let $G$ be the Galois group of $\Q(\zeta)/\Q$. Let $d\in\N,\quad p-1\equiv 0\modu d$. Let 
$G_d$ be the subgroup of the cyclic group $G$. Then $G_d$ is  of order $\frac{p-1}{d}$.
If $\sigma$ generates $G$, then $\sigma^d$ generates $G_d$.
\item
Let $\mathbf b$ be an ideal of $\Z[\zeta]$, not principal and with $\mathbf b^p$ principal.
\item
Let  $c_i= Cl(\sigma^i(\mathbf b)),\quad i=0,\dots,p-2,$  be the class of $\sigma^i(\mathbf b)$ in the $p$-class group of $\Q(\zeta)$.
\end{itemize}
\begin{itemize}
\item 
Recall that $Cl(\mathbf b)$ is the class of the ideal $\mathbf b$ of $\Z[\zeta]$.
Observe that exponential notations $\mathbf b^\sigma$ can be used indifferently in the sequel.
With this notation, we have
\begin{itemize}
\item
$\mathbf b^{\sigma^d}=\sigma^d(\mathbf b)$.
\item
For $\lambda\in {\bf F}_p$, we have 
$ \mathbf b^{\sigma+\lambda} = \mathbf b^\lambda\times  \sigma(\mathbf b).$
\item
Let $P(\sigma)= 
\sigma^m+\lambda_{m-1}\sigma^{m-1}+\dots+\lambda_1\sigma+\lambda_0
 \in {\bf F}_p[\sigma]$; then 
$ \mathbf b^{P(\sigma)} =
\sigma^m(\mathbf b)\times \sigma^{m-1}(\mathbf b)
^{\lambda_{m-1}}\times\dots\times\sigma(\mathbf b)^{\lambda_1}
\times \mathbf b^{\lambda_0}$.
\item
Let us note $ \mathbf b^{P(\sigma)} \simeq \Z[\zeta]$,  if the ideal 
$\sigma^m(\mathbf b)\times \sigma^{m-1}(\mathbf b)^{\lambda_{m-1}}\dots \sigma(\mathbf b)^{\lambda_1}\times \mathbf b^{\lambda_0}$ is principal.
\item
Let $P(\sigma), Q(\sigma)\in {\bf F}_p[\sigma]$;
if $\mathbf b^{P(\sigma)}\simeq \Z[\zeta]$, 
then $ \mathbf b^{Q(\sigma)\times P(\sigma)}\simeq \Z[\zeta]$.
\item
Observe that trivially $ \mathbf b^{\sigma^{p-1}-1} \simeq \Z[\zeta].$
\end{itemize}
\item 
There exists a  monic minimal  polynomial $P_{r_d}(V)\in {\bf F}_p[V]$, polynomial ring of the indeterminate $V$  verifying the relation, for $V=\sigma^d$: 
\begin{equation}\label{e30101}
 \mathbf b^{P_{r_d}(\sigma^d)} 
\simeq \Z[\zeta].
\end{equation}
This  minimality implies that,   for all polynomials $R(V)\in {\bf F}_p(V), \quad R(V)\not=0,\quad 
deg(R(V))<deg(P_{r_d}(V))$, we have $  \mathbf b^{ R(\sigma^d)}\not\simeq \Z[\zeta]$.
It means, with an other formulation in term of ideals, that 
$\prod_{i=0}^{r_d} \sigma^{i d}(\mathbf b)^{\lambda_{i,d}}$ is a principal ideal and 
that $\prod_{i=0}^{\alpha} \sigma^{i d}(\mathbf b)^{\beta_i}$ is not principal when $\alpha<r_d$ and $\beta_i,\quad i=0,\dots,\alpha$, are not all simultaneously null.
\item
$P_{r_d}(U)$ is the {\bf minimal polynomial} of the indeterminate $U$ with  
$P_{r_d}(\sigma^d)\in{\bf F}_p[G_d]$ {\bf annihilating} the ideal class of 
$\mathbf b$.
\end{itemize}
%
\subsection{Representations of Galois group $Gal(\Q(\zeta)/\Q)$ in characteristic $p$.}\label{s20112}
\label{s07061a}

In this subsection we give some general properties of representations of $G=Gal(\Q(\zeta)/\Q)$ in characteristic $p$ and obtain some results on the structure of  the $p$-class group of $\Q(\zeta)$.
Observe that we never use characters theory.
%
\begin{lem} \label{l04101}
Let $d\in\N, \quad p-1\equiv 0\modu d$.
Let $V$ be an indeterminate. Then the minimal  polynomial $P_{r_d}(V)$ 
with $P_{r_d}(\sigma^d)\in{\bf F}_p[G_d]$  annihilating ideal class of $\mathbf b$ 
verifies the factorization
\begin{displaymath}
P_{r_d}(V)=\prod_{i=1}^{r_d} (V-\mu_{i, d}),\quad \mu_{i,d}\in{\bf F}_p,\quad
 i_1\not= i_2 \Rightarrow \mu_{i_1}\not= \mu_{i_2}.
\end{displaymath}
\begin{proof}
Let us consider the polynomials $A(V)=V^{p-1}-1$ and $P_{r_d}(V)\in {\bf F}_p[V]$. It is possible to divide the polynomial  $A(V)$ by $P_{r_d}(V)$ in the polynomial ring 
${\bf F}_p[V]$ to obtain
\begin{displaymath}
\begin{split}
& A(V)=P_{r_d}(V)\times Q(V)+R(V),\quad Q(V),R(V)\in{\bf F}_p[V],\\ 
& d_R=deg_V(R(V))<  r_d = deg_V(P_{r_d}(V)).
\end{split}
\end{displaymath}
For $V=\sigma^d$, we get 
 $ \mathbf b^{\sigma^{d(p-1)}-1} \simeq \Z[\zeta]$ and 
$  \mathbf b^{P_{r_d}(\sigma^d)} \simeq \Z[\zeta]$, so  $  \s^{R(\sigma^d)} \simeq \Z[\zeta]$.
Suppose that $R(V)=\sum_{i=0}^{d_R} R_i V^i,\quad R_i\in {\bf F}_p$, is not identically null; then, it leads to the relation
\begin{displaymath}
\sum_{i=0}^{d_R} R_i \sigma ^{d i}(\mathbf b)\simeq \Z[\zeta],
\end{displaymath}
where the $R_i$ are not all zero,
with $d_R <r_d$, which contradicts the minimality of the polynomial $P_{r_d}(V)$.
Therefore, $R(V)$ is identically null and we have
\begin{displaymath}
V^{p-1}-1= P_{r_d}(V)\times Q(V).
\end{displaymath}
The factorization of $V^{p-1}-1$ in ${\bf F}_p[V]$ is 
$V^{p-1}-1
=\prod_{i=1}^{p-1} (V-i)$.
The factorization is unique in the euclidean ring ${\bf F}_p[V]$ and so  
$P_{r_d}(V)=\prod_{i=1}^{r_d} (V-\mu_{i,d}),\quad \mu_{i,d}\in{\bf F}_p,
\quad  i_1\not= i_2 \Rightarrow \mu_{i_1}\not= \mu_{i_2}$, which achieves the proof.
\end{proof}
\end{lem}
%
***************************************

\begin{lem}\label{l22101}
Let $d\in\N,\quad p-1\equiv 0\modu d$.
Let $U, W$ be two indeterminates. Let $P_{r_1}(U)$ be the minimal polynomial with $P_{r_1}(\sigma)
\in {\bf  F}_p[G]$ annihilating the ideal class of $\mathbf b$.
Let $P_{r_d}(W)$ be the minimal polynomial with $P_{r_d}(\sigma^d)
\in {\bf  F}_p[G_d]$  annihilating the ideal class of $\mathbf b$.
Then  
\begin{enumerate}
\item
$ P_{r_1}(U)=\prod_{i=1}^{r_1} (U-\mu_i),\quad \mu_i\in{\bf F}_p.$
\item
$ P_{r_d}(U^d)=\prod_{i=1}^{r_d} (U^d-\mu_i^d)
= P_{r_1}(U)\times Q_d(U),\quad r_d\leq r_1,
\quad Q_d(U)\in{\bf F}_p[U]$.
\item
The $p$-ranks $r_1$ and $r_d$ verify the inequalities
\begin{equation}\label{e011181}
r_d\times d\geq  r_1\geq r_d.
\end{equation}
\item
Let $K_d$ be the intermediate field $\Q\subset K_d\subset\Q(\zeta),
\quad [K_d:\Q]=d$. Suppose that $p$ does not divide the class number of $K_d/\Q$;
then $\mu_i^d\not=1$ for $i=1,\dots,r_d$. In particular   $\mu_i\not=1$ for $i=1,\dots,r_1$.
\end{enumerate}
\begin{proof}$ $
\begin{itemize}
\item
Observe, at first, that $deg_U(P_{r_d}(U^d))=d\times r_d\geq r_1$: if not, for the polynomial $P_{r_d}(U^d)$ seen in the indeterminate $U$, whe should have
$deg_U(P_{r_d}(U^d))<r_1$ and 
$P_{r_d}(\sigma^d) \circ \mathbf b \simeq \Z[\zeta]$ and, as previously, the polynomial  $P_{r_d}(U^d)$ of the indeterminate $U$ should be identically null.
\item 
We  apply euclidean algorithm in the polynomial ring ${\bf F}_p[U]$ of the indeterminate $U$.
Therefore,  
\begin{displaymath}
\begin{split}
& P_{r_d} (U^d)=P_{r_1}(U)\times Q(U)+R(U),\quad Q(U),R(U)\in{\bf F}_p[U],\\
& deg(R(U))<deg(P_{r_1}(U)).
\end{split}
\end{displaymath}
But we have $\mathbf b^{P_{r_d}(\sigma^d)} \simeq\Z[\zeta],
\quad  \mathbf b^{P_{r_1}(\sigma )} \simeq \Z[\zeta]$, therefore 
$\mathbf b^{R( \sigma)}\simeq \Z[\zeta]$.
Then, similarly to proof of lemma \ref{l04101} p.\pageref{l04101}, $R(U)$ is identically null and 
$P_{r_d}(U^d)= P_{r_1}(U) \times Q(U)$.
\item
Applying lemma \ref{l04101} p.\pageref{l04101},  we obtain
\begin{displaymath}
\begin{split}
&P_{r_1}(U)=\prod_{i=1}^{r_1} (U-\mu_{i}),\quad \mu_{i}\in{\bf F}_p,\\
&P_{r_d}(U^d)=\prod_{i=1}^{r_d} (U^d-\mu_{i,d}),\quad \mu_{i,d}\in{\bf F}_p.
\end{split}
\end{displaymath}
\item
Then, we get 
\begin{displaymath}
P_{r_d}(U^d)=\prod_{i=1}^{r_d} (U^d-\mu_{i,d})=\prod_{i=1}^{r_1} (U-\mu_i)\times Q(U).
\end{displaymath}
There exists at least one $i, \quad 1\leq i \leq r_d$, such that 
$(U^d-\mu_{i,d})=(U-\mu_1)\times Q_1(U)$: 
if not, for all $i=1,\dots,r_d$, we should have
$U^d- \mu_{i,d}\equiv R_i\modu (U-\mu_1),\quad R_i\in {\bf F}_p^*$, a contradiction because $\prod_{i=1}^{r_d} R_i\not=0$. 
We have $\mu_{i,d}=\mu_1^d$: if not $U-\mu_1$ should divide $U^d-\mu_{i,d}$ and $U^d-\mu_1^d$ and also 
$U-\mu_1$ should divide 
$(\mu_{i,d}-\mu_1^d)\in {\bf F}_p^*$, a contradiction.
Therefore, there exists at least one $i, \quad 1 \leq i\leq r_d$, such that $\mu_{i,d}=\mu_1^d$ and $U^d-\mu_{i,d}= U^d-\mu_1^d= (U-\mu_1)\times Q_1(U)$. 

Then, generalizing to $\mu_{i,d}$ for all $i=1,\dots,r_d$, we get with a certain reordering of index $i$ 
\begin{displaymath}
P_{r_d}(U^d)=\prod_{i=1}^{r_d} (U^d-\mu_i^d)
=\prod_{i=1}^{r_1} (U-\mu_i)\times Q(U).
\end{displaymath}
\item
We have 
\begin{displaymath}
P_{r_d}(U^d)=\prod_{i=1}^{r_d} (U^d-\mu_i^d).
\end{displaymath}
This relation leads to
\begin{displaymath}
P_{r_d}(U^d)=\prod_{i=1}^{r_d}\prod_{j=1}^d(U-\mu_i\mu_d^j),
\end{displaymath}
where $\mu_d\in {\bf F}_p,\quad \mu_d^d=1$.
We have shown that $P_{r_d}(U^d)=P_{r_1}(U)\times Q_d(U)$ and so
$deg_U(P_{r_d}(U))=d\times r_d\geq r_1$;  thus $d\times r_d\geq r_1$.
\item
We finish by the proof of item 4): suppose that, for some $i,\quad 1\leq i\leq r_d$, we have $\mu_i^d=1$ and search for a contradiction:
there exists, for the indeterminate V,  a polynomial $P_1(V)\in {\bf F}_p(V)$ such that 
$P_{r_d}(V)=(V-\mu_i^d)\times P_1(V)=(V-1)\times P_1(V)$. 
But for $V= \sigma^d$, we have $ \mathbf b ^{P_{r_d}(\sigma^d) }\simeq \Z[\zeta]$,
so $ \mathbf b ^{(\sigma^d P_{1}(\sigma^d)- P_1(\sigma^d))}\simeq\Z[\zeta]$. So, $\mathbf b^{P_1(\sigma^d)} $ is the class of an  ideal $\mathbf c$ of $\Z[\zeta]$ with $Cl(\sigma^d(\mathbf c))=Cl(\mathbf c)$;
then  $Cl(\sigma^{2d}(\mathbf c))
=Cl(\sigma^d(\mathbf c))=Cl(\mathbf c)$.
Then $Cl(\sigma^d(\mathbf c)\times\sigma^{2d}(\mathbf c)\times\dots
\times\sigma^{(p-1)d/d}(\mathbf c))
=Cl(\mathbf c^{(p-1)/d})$.
Let $\tau=\sigma^{d}$; then
$Cl(\tau(\mathbf c)\times\tau^{2 }(\mathbf c)\times\dots\times
\tau^{(p-1)/d}(\mathbf c))=Cl(\mathbf c^{(p-1)/d})$;
Then we deduce that 
$Cl(N_{\Q(\zeta)/K_d} (\mathbf c))=Cl(\mathbf c^{(p-1)/d})$ and thus 
$\mathbf c$ is a principal ideal because the ideal  $N_{\Q(\zeta)/K_d} (\mathbf c)$
of $K_d$ is principal,  (recall that, from hypothesis, $p$ does not divide $h(K_d/\Q)$); so $ \mathbf b^{P_1(\sigma^d)} \simeq\Z[\zeta]$, which contradicts the minimality of 
the minimal polynomial equation $ \mathbf b ^{P_{r_d}(\sigma^d)}\simeq\Z[\zeta]$ because, for the indeterminate $V$, we would have  $deg(P_{1}(V))<deg(P_{r_d}(V))$, which achieves the proof.
\end{itemize}
\end{proof}
\end{lem}
{\bf Remark:} As an example, the item 4) says that:
\begin{enumerate}
\item
If $d=2$, then classically $p\ \not|\ h(K_2/\Q)$ and so item 4) shows that $\mu_i\not=-1$: there is no ideal 
$\mathbf b$ whose class belongs to $C_p$ which is  annihilated by $\sigma-u_{(p-1)/2}=\sigma+1$.
\item if $h^+\not\equiv 0\modu p$, (Vandiver's conjecture) then 
$\mu_i^{(p-1)/2}=-1$ for $i=1,\dots,r_1$.
\end{enumerate}
%
We summarize  results obtained in: 
\begin{lem}\label{t31101a} 
Let $\mathbf b$ be an ideal of $\Z[\zeta],\quad \mathbf b^p\simeq \Z[\zeta],\quad 
\mathbf b\not \simeq \Z[\zeta]$. 
Let $d\in\N,\quad p-1\equiv 0\modu d$.
Let $U, W$ be two indeterminates. Let $P_{r_1}(U)$ be the minimal polynomial with $P_{r_1}(\sigma)
\in {\bf  F}_p[G]$ annihilating the ideal class of $\mathbf b$.
Let $P_{r_d}(W)$ be the minimal polynomial with $P_{r_d}(\sigma^d)
\in {\bf  F}_p[G_d]$  annihilating the ideal class of $\mathbf b$.
Then there exists $\mu_1,\mu_2,\dots,\mu_{r_1}\in{\bf F}_p$,
with $i\not=i^\prime\Rightarrow \mu_i\not=\mu_i^\prime$,  such that, for the indeterminate $U$,
\begin{itemize}
\item
the minimal polynomials $P_{r_1}(U)$ and $P_{r_d}(U^d)$ are respectively given by
\begin{displaymath}
\begin{split}
&P_{r_1}(U)=\prod_{i=1}^{r_1} (U-\mu_i),\\
&P_{r_d}(U^d)=\prod_{i=1}^{r_d} (U^d-\mu_i^d),\quad r_d\leq r_1,\\
&P_{r_1}(U) \ | \ P_{r_d}(U^d).
\end{split}
\end{displaymath}
\item
The coefficients of $P_{r_d}(U^d)$ are explicitly computable by
\begin{displaymath}
\begin{split}
& P_{r_d}(U^d)=\\
& U^{d r_d} -S_1(d)\times U^{d(r_d-1)}+S_2(d)\times U^{d(r_d-2)}
+\dots+(-1)^{{r_d}-1} S_{{r_d}-1}(d)\times U^d
+(-1)^{r_d} S_{r_d}(d),\\
& S_0(d)=1,\\
& S_1(d)=\sum_{i=1,\dots, r_d} \mu_i^d,\\
& S_2(d)=
\sum_{1\leq i_1<i_2\leq r_d} \mu_{i_1}^d\mu_{i_2}^d,\\
&\vdots\\
&S_{r_d}(d)=\mu_1^d\mu_2^d\dots\mu_{r_d}^d.
\end{split}
\end{displaymath}
\item
Then the ideal 
\begin{equation}\label{e04112}
\prod_{i=0}^{r_d} \sigma^{di}(\mathbf b)^{(-1)^{r_d-i}\times S_{r_d-i}(d)}
 =\mathbf b^{P_{r_d}(\sigma^d)}
\end{equation} 
is a principal ideal.
\end{itemize}
\end{lem}
\paragraph{Remark:} For other annihilation methods of $Cl(\Q(\zeta)/\Q)$ more involved, see for instance Kummer, in Ribenboim \cite{rib} p 119, (2C) and (2D) and Stickelberger in Washington \cite{was} p 94 and 332.
%
\subsection{On the structure of the $p$-class group of subfields of $\Q(\zeta)$}
\label{s108031}
In  this subsection we get  several  results on  the structure of the $p$-class group  of $\Q(\zeta)$ and on class number $h$ of $\Q(\zeta)$:
\begin{itemize}
\item
A formulation, with our notations,  of a Ribet's  result on irregularity index.
\item
Let $d,g\in\N$ coprime with $d\times g=p-1$. 
For groups generated by the action of  Galois groups $G$ and of subgroups 
$G_d, G_g$ of $G$  on ideals $\mathbf b$ of $\Q(\zeta)$,
an inequality between  degrees $r_1, r_d, r_g$  in the indeterminate $X$  of minimal polynomials 
$P_{r_1}(X), P_{r_d}(X), 
P_{r_g}(X)\in {\bf F}_p[X]$, with $P_{r_1}(\sigma), P_{r_d}(\sigma^d), P_{r_g}(\sigma^g)$  annihilating ideal class of $\mathbf b$.
\item
Some $\pi$-adic congruences  connected to structure of $p$-class group $C_p$ of $\Q(\zeta)$.
\end{itemize}
%
\subsubsection{Some definitions and notations}
\begin{itemize}
\item
Recall that:  
\begin{itemize}
\item
$r_p$ is the $p$-rank of the class group of $\Q(\zeta)$.
\item
$C_p$ is the subgroup of   exponent $p$ of the $p$-class group of $\Q(\zeta)$.
\end{itemize}
\item
The $\Q$-isomorphism $\sigma$ of $\Q(\zeta)$ generates  $G=Gal(\Q(\zeta)/\Q)$,  Galois group of the field $\Q(\zeta)$. 
For $d\ |\ p-1$, let  $G_d$ be the subgroup of $d$ powers $\sigma^{di}$ of elements $\sigma^i$ of $G$.
This group is  of order $\frac{p-1}{d}$.
\item
Suppose that $r_p>0$.
There exists an ideal class with  representants $\mathbf b\subset\Z[\zeta]$, 
with $\mathbf b^p\simeq \Z[\zeta],\quad \mathbf b\not\simeq\Z[\zeta]$, which verifies, in term of representations, for some  ideals  $\mathbf b_i$ of $\Z[\zeta],\quad i=1,\dots,r_p$, 
\begin{equation}\label{e107162}
\begin{split}
& \mathbf b\simeq\prod_{i=1}^{r_p} \mathbf  b_i,\\
& \mathbf b_i^p\simeq \Z[\zeta],\quad \mathbf b_i\not\simeq \Z[\zeta],\quad i=1,\dots,r_p,\\
& \sigma(\mathbf b_i)\simeq\mathbf b_i^{\mu_i},
\quad \mu_i\in{\bf F}_p,\quad \mathbf b_i+\pi=\Z[\zeta],\quad i=1,\dots,r_p,\\
& C_p= \oplus_{i=1}^{r_p} <Cl(\mathbf b_i)>,\\
& P_{r_1}(U)=\prod_{i=1}^{r_1} (U-\mu_i),\quad \mathbf b^{P_{r_1}(\sigma)}
\simeq \Z[\zeta],\quad 1\leq r_1\leq r_p, 
\end{split}
\end{equation}
where $P_{r_1}(U)$ is the minimal polynomial in the indeterminate $U$ for the action of $G$ on the ideal  $\mathbf b$, such that $P_{r_1}(\sigma)\in{\bf F}_p[G]$ annihilates the ideal class of $\mathbf b$, 
see theorem \ref{t31101a} p \pageref{t31101a}.
Recall that it is possible to  encounter the case  $\mu_i=\mu_j$ in the set
$\{\mu_1,\dots,\mu_{r_p}\}$; by opposite if $U-\mu_i$ and $U-\mu_j$ divide the minimal polynomial $P_{r_1}(U)$ then $\mu_i\not=\mu_j$. Therefore $r_1$ is  the degree of the minimal polynomial $P_{r_1}(U)$.
\item
With a certain indexing assumed in the sequel, the ideals classes  $Cl(\mathbf b_i)\in C_p^-$ for $i=1,\dots,r_p^-$,
and ideal classes  $Cl(\mathbf b_i)\in C_p^+$ for $i=r_p^-+1,\dots, r_p$. 
\begin{itemize}
\item
The ideal $\mathbf b$ verifies $\mathbf b\simeq \mathbf b^-\times\mathbf b^+$ where $\mathbf b^-$ and 
$\mathbf b^+$ are two ideals of $\Q(\zeta)$ with  
$Cl(\mathbf b^-)\in C_p^-$ and $Cl(\mathbf b^+)\in C_p^+$.
\item
With this notation, the minimal polynomial $P_{r_1}(U)$ factorize in a factor corresponding to $C_p^-$ and a factor corresponding to $C_p^+$, with:
\begin{equation}\label{e301021}
P_{r_1}(U)=P_{r_1^-}(U)\times P_{r_1^+}(U),\quad r_1=r_1^-+r_1^+.
\end{equation}
\item
$P_{r_1^-}(U)$ is the minimal polynomial with $P_{r_1^-}(\sigma)\in{\bf F}_p[G]$ annihilating the class of ideal $\mathbf b^-\in C_p^-$.
\item
$P_{r_1^+}(U)$ is the minimal polynomial with $P_{r_1^+}(\sigma)\in{ \bf F}_p[G]$ annihilating the class of ideal $\mathbf b^+\in C_p^+$.
\end{itemize}
\item
Let us denote $M_{r_1}=\{\mu_i\quad|\quad i=1,\dots,r_1\}$.
\item 
Let $d\in\N,\quad d\ |\ p-1,\quad 2\leq d\leq\frac{p-1}{2}$. Let $K_d$ be the field
$\Q\subset K_d\subset\Q(\zeta),\quad [K_d:\Q]=d$.
\item 
Let $P_{r_d}(V)$ be the minimal polynomial in the indeterminate $V$ of the action of the group $G_d$ on the ideal class group $<\mathbf b>$ of order $p$,  such that $P_{r_d}(\sigma^d)\in{\bf F}_p[G_d]$ annihilates ideal class of $\mathbf b$.
Let $r_d$ be the degree of $P_{r_d}(V)$. 
\end{itemize}
%
\subsubsection{ On the irregularity index}
Recall that $r_p$ is the $p$-rank of the group $C_p$.
The irregularity index is the number 
\begin{displaymath}
i_p=Card\{B_{p-1-2m}\ | \ B_{p-1-2m}\equiv 0\modu p,\quad 1\leq m\leq\frac{p-3}{2}\},
\end{displaymath}
where $B_{p-1-2m}$ are even Bernoulli Numbers.
The next theorem connects irregularity index and degree  $r_1^-$ of minimal polynomial $P_{r_1^-}(U)$ defined in relations (\ref{e107162}) p. \pageref{e107162} and \ref{e301021} p. \pageref{e301021}.
%
\begin{thm}{ *** }\label{t204251}
With meaning of degree $r_1^-$ of minimal polynomial  $P_{r_1^-}(U)$ defined in relation (\ref{e301021}) p. \pageref{e301021}, 
then the irregularity index is equal to the degree $r_1^-$ and   verifies :
\begin{equation}\label{e301041}
r_p^--r_p^+\leq i_p=r_1^-\leq r_p^-. 
\end{equation}
\begin{proof}
Let us consider in relation (\ref{e107162}) the set of ideals  $\{\mathbf b_i\ |\ i=1,\dots,r_p\}$.
The result of Ribet using theory of modular forms \cite{rie} mentionned in Ribenboim \cite{rib} (8C) p 190 can be formulated, with our notations, 
\begin{equation}\label{e204251}
B_{p-1-2m}\equiv 0\modu p \Leftrightarrow \exists i,\quad 1 \leq i\leq r_p,\quad 
\mathbf b_i^{\sigma-u_{2m+1}}\simeq \Z[\zeta].
\end{equation}
There exists at least one such $i$, but it is possible for $i\not=i^\prime$ that 
$\mathbf b_i^{\sigma-u_{2m+1}}\simeq \mathbf b_{i^\prime}^{\sigma-u_{2m+1}}\simeq \Z[\zeta]$.
\begin{itemize}
\item 
The relation (\ref{e204251}) p. \pageref{e204251} implies that $i_p=r_1^-$.
\item
The inequality (\ref{e301041}) p. \pageref{e301041} is an immediate consequence of independant forward structure 
theorem \ref{t203021} p. \pageref{t203021}.
\end{itemize}
\end{proof}
\end{thm}
%
%

%
\subsubsection{Inequalities involving  degrees $r_1, r_d, r_g $ of minimal polynomials $P_{r_1}(V),P_{r_d}(V),P_{r_g}(V)$ annihilating ideal $\mathbf b$.}
In this subsection, we always assume that $\mathbf b$ is defined by relation  (\ref{e107162}) p. \pageref{e107162}.

Let $p$ be an odd prime.
Let $d, g\in \N$, with $gcd(d,g)=1$ and $ d\times g = p-1$.
Recall that    $r_1$, $r_d$ and $r_g$ are the degrees of the minimal polynomials
$P_{r_1}(V), P_{r_d}(V), P_{r_g}(V)$ of the indeterminate $V$ with
$\mathbf b^{P_{r_1}(\sigma)}\simeq \mathbf b^{P_{r_d}(\sigma^d)}
\simeq \mathbf b^{P_{r_g}(\sigma^g)}\simeq \Z[\zeta]$. 
%
The next theorem is a relation between the three degree  $r_1, r_d$ and $r_g$. 
\begin{thm}{ *** }\label{l108081}
Let $d, g\in\N, \quad gcd(d,g)=1,\quad d\times g = p-1$. Suppose that $r_d\geq 1$ and $r_g\geq 1$. 
Then
\begin{equation}\label{e108081}
r_d\times r_g\geq r_1.
\end{equation}
and if $r_d=1$ then $r_g=r_1$.
\begin{proof} $ $
\begin{itemize}
\item
Let us consider the minimal polynomials $P_{r_d}(U^d)=\prod_{i=1}^{r_d}
(U^d-\mu_i^d)$ and $P_{r_g}(U^g)=\prod_{i=1}^{r_g}(U^g-\nu_j^g)$ of the indeterminate $U$ with $\mathbf b^{P_{r_d}(\sigma^d)}\simeq \Z[\zeta]$ and
$\mathbf b^{P_{r_g}(\sigma^g)}\simeq \Z[\zeta]$.
\item
From lemma \ref{l22101} p.\pageref{l22101}, 
we have seen that $P_{r_1}(U)\ |\ P_{r_d}(U^d)$ and that similarly
$P_{r_1}(U)\ |\ P_{r_g}(U^g)$, thus $P_{r_1}(U)\ |\ gcd(P_{r_d}(U^d), P_{r_g}(U^g))$.
\item
Let $M_{r_1}=\{\mu_i\quad|\quad i=1,\dots,r_1\}$.
Let us define the sets
\begin{displaymath}
C_1(\mu_i)=\{\mu_i\times \alpha_j\quad|\quad\alpha_j^d=1,\quad j=1,\dots,d\}\cap M_{r_1}, \quad i=1,\dots,r_d.
\end{displaymath}
Let us define in the same way the sets
\begin{displaymath}
C_2(\nu_i)=\{\nu_i\times \beta_j\quad|\quad\beta_j^g=1,\quad j=1,\dots,g\}\cap M_{r_1}, \quad i=1,\dots,r_g.
\end{displaymath}
\item
We have proved in lemma \ref{l22101} p.\pageref{l22101} that $P_{r_1}(U)\ |\ P_{r_d}(U^d)$. Therefore the sets  
$C_1(\mu_i),\quad i=1,\dots,r_d$, are a partition of $M_{r_1}$ and 
$r_1=\sum_{i=1}^{r_d} Card(C_1(\mu_i))$.
\item
In the same way $P_{r_1}(U)\ |\ P_{r_g}(U^g)$. Therefore the sets  
$C_2(\nu_i),\quad i=1,\dots,r_g$, are a partition of $M_{r_1}$ and 
$r_1=\sum_{i=1}^{r_g} Card(C_2(\nu_i))$.
\item
There exists at least one $i\in\N,\quad 1\leq i\leq r_d$, such that 
$Card(C_1(\mu_i))\geq \frac{r_1}{r_d}$. For this $i$, let 
$\nu_1 =\mu_i\times \alpha_1,\quad \alpha_1^d=1,\quad \nu_1\in M_{r_1}$ and, in the same way, let 
$\nu_2 =\mu_i\times \alpha_2,\quad \alpha_2^d=1,\quad \nu_2\in M_{r_1},\quad \nu_2\not=\nu_1$.
We have $\nu_1^g\not=\nu_2^g$ : if not we should simultaneously have
$\alpha_1^d=\alpha_2^d$ and $\alpha_1^g=\alpha_2^g$, which should imply, from
$gcd(d,g)=1$, that $\alpha_1=\alpha_2$, contradicting $\nu_1\not=\nu_2$ and therefore we get $C_2(\nu_1)\not=C_2(\nu_2)$.
\item
Therefore, extending the same reasoning to all elements of $C_1(\mu_i)$, we get  $\frac{r_1}{r_d}\leq Card(C_1(\mu_i))\leq r_g$, which leads to the result.
\item
If $r_d=1$ then $r_g\geq r_1$ and in an other part $r_g\leq r_1$ and so $r_g=r_1$.
\end{itemize}
\end{proof}
\end{thm}
%
\paragraph{Remarks:} 
\begin{itemize}
\item
As an example, consider an odd prime $p$ verifying $p\not\equiv 1\modu 4$. Suppose also that $h^+\not\equiv 0\modu p$. Then $P_{r_{(p-1)/2}}(\sigma)=\sigma^{(p-1)/2}+1=U+1$ for the indeterminate $U=\sigma^{(p-1)/2}$.
Therefore $r_{(p-1)/2}=1$ and thus $r_2=r_1$.
\item
Observe that $1\leq r_d<r_1$ implies that $r_g>1$.
\end{itemize}
%
%
%
\subsubsection {On Stickelberger's ideal in field $\Q(\zeta)$}
In this subsection, we give a result   resting on the annihilation of class group of $\Q(\zeta)$ by Stickelberger's ideal.
\begin{itemize}
\item 
Let us denote $\mathbf a\simeq \mathbf c$ when the two ideals $\mathbf a$ and $\mathbf c$ of $\Q(\zeta)$ are in the same ideal class.
\item
Let $G=Gal(\Q(\zeta)/\Q)$.
\item
Let $\tau_a :\zeta \rightarrow \zeta^a,\quad a=1,\dots,p-1$, be the $p-1$ $\Q$-isomorphisms of
the field $\Q(\zeta)/\Q$.
\item
Recall that $u$ is a primitive root $\modu p$,  and that 
$\sigma : \zeta \rightarrow \zeta^u$ 
is a $\Q$-isomorphism of the field $\Q(\zeta)$ 
which generates $G$.  Recall that, for $i\in \N$, then we denote $u_i$ for 
$u^i \modu p$ and $1\leq u_i\leq p-1$.
\item
Let $\mathbf b$ be the not principal ideal defined in relation (\ref{e107162}) p.\pageref{e107162}.
Let $P_{r_1}(\sigma)\in {\bf F}_p[G]$ be the polynomial of minimal degree such that
$P_{r_1}(\sigma)$ annihilates $\mathbf b$, so such that
$\mathbf b^{P_{r_1}(\sigma)}$ is principal ideal, see lemma \ref{l04101} p.\pageref{l04101}, and so
\begin{displaymath}
P_{r_1}(\sigma)=\prod_{i=1}^{r_1} (\sigma-\mu_i),\quad \mu_i\in {\bf F}_p,
\quad i\not=i^\prime \Rightarrow \mu_i\not=\mu_i^\prime.
\end{displaymath}
\end{itemize}
In the next result we shall explicitly use the annihilation of class group of $\Q(\zeta)$ by the Stickelberger's ideal.
%
\begin{lem}\label{l103021}
Let $P_{r_1}(U)=\prod_{i=1}^{r_1}(U-\mu_i)$ be the polynomial of the indeterminate $U$, of minimal degree,  such that $\mathbf b^{P_{r_1}(\sigma)}$ is principal.  Then $\mu_i\not= u,\quad i=1,\dots,r_1$.
\begin{proof}$ $
\begin{itemize}
\item
Let $i\in\N, \quad 1\leq i\leq r_1$.
From relation (\ref{e107162}) p.\pageref{e107162}, there exists   ideals $\mathbf b_i\in \Z[\zeta],\quad i=1,\dots,r_p$, not principal and such that $\mathbf b=\prod_{i=1}^{r_p} \mathbf b_i$,  with $\mathbf b_i^{\sigma-\mu_i}$  principal.
\item
Suppose that $\mu_i=u$, and search for a contradiction:
Let us consider  $\theta=\sum_{a=1}^{p-1}\frac{a}{p}\times\tau_a^{-1}\in\Q[G]$. 
Then $p\theta\in\Z[G]$ and the ideal
$\mathbf b^{p\theta}$ is principal from Stickelberger's theorem, see for instance Washington \cite{was}, theorem 6.10 p 94.
\item
We can set  $a=u^{m},\quad a=1,\dots,p-1$, and $m$ going through all the set $\{0,1,\dots,p-2\}$, because $u$ is a primitive root $\modu p$. 
Then $\tau_a :\zeta\rightarrow \zeta^a$ and so 
$\tau_a^{-1} :\zeta \rightarrow \zeta^{(a^{-1})}=\zeta^{((u^m)^{-1})}
=\zeta^{(u^{-m})}
=\zeta^{(u^{p-1-m})}=\sigma^{p-1-m}=\sigma^{-m}$.
\item
Therefore,  
$p\theta=\sum_{m=0}^{p-2} u^m\sigma^{-m}$.
The element $\sigma-\mu_i=\sigma-u$ annihilates the class of $\mathbf b_i$ and also the element 
$u\times\sigma^{-1}-1$ annihilates the class of $\mathbf b_i$. Therefore 
$u^{m}\sigma^{-m}-1,\quad m=0,\dots,p-2,$ annihilates the class of $\mathbf b_i$ and finally $p-1$ annihilates the class of $\mathbf b_i$, so $\mathbf b_i^{p-1}$ is principal, but $\mathbf b_i^p$ is also principal, and finally $\mathbf b_i$ is principal which contradicts our hypothesis and achieves the proof.
\end{itemize}
\end{proof}
\end{lem}
%
%
%
\subsection{$\pi$-adic congruences  connected  to $p$-class group $C_p$}\label{s108311}
In a first subsection, we examine  the case of relative $p$-class group $C_p^-$.
In a second  subsection, we examine   the case of $p$-class group $C_p^+$.
In last subsection, we summarize our results to all $p$-class group $C_p$.
These important congruences (subjective) characterize structure of $p$-class group.
%
\subsubsection{$\pi$-adic congruences  connected  to relative $p$-class group $C_p^-$}
In this subsection , we shall describe  some $\pi$-adic congruences  connected  to $p$-relative class group $C_p^-$.
\paragraph{Some definitions and a preliminary result}
\begin{itemize}
\item
Let $C_p$ be the subgroup of exponent $p$ of the $p$-class group of $\Q(\zeta)$.
\item
Let $r_p$ be the $p$-rank of $C_p$, let  $r_p^+$ be the $p$-rank of $C_p^+$ and $r_p^-$ be the relative $p$-rank of $C_p^-$.
Let us recall the structure of  the ideal $\mathbf B$ already defined in relation (\ref{e107162}) p.\pageref{e107162}:
\begin{equation}\label{e203031}
\begin{split}
& \mathbf B=\mathbf b_1\times\dots\times\mathbf b_{r_p^-}
\times \mathbf b_{r_p^-+1}\times\dots\times \mathbf b_{r_p},\\
& C_p=\oplus_{i=1}^{r_p} <Cl(\mathbf b_i)>,\\
& \mathbf b_i^p\simeq\Z[\zeta],\quad \mathbf b_i\not\simeq \Z[\zeta],
\quad i=1,\dots,r_p,\\
&\sigma(\mathbf b_i)\simeq \mathbf b_i^{\mu_i},\quad\mu_i\in{\bf F}_p^*,
\quad i=1,\dots,r_p,\\
& Cl(\mathbf b_i)\in C_p^-,\quad i=1,\dots,r_p^-,\\
& Cl(\mathbf b_i)\in C_p^+,\quad i=r_p^-+1,\dots,r_p,\\
& \mathbf B^{P_{r_1}(\sigma)}\simeq\Z[\zeta],\\
& (\frac{\mathbf B}{\overline{\mathbf B}})^{P_{r_1^-}(\sigma)}\simeq\Z[\zeta].
\end{split}
\end{equation}
(Observe that we replace here notation $\mathbf b$ by $\mathbf B$ to avoid conflict of notation in the sequel.)
In the sequel, we are using also the natural integers $m_i$, with $0\leq m_i\leq \frac{p-3}{2}$, defined by $\mu_i=u_{2m_i+1}=u^{2m_i+1}\modu p$.
\item
Recall that it is possible to have $\mu_i=\mu_j=\mu$: observe that, in that case,  the decomposition
$<Cl(\mathbf b_i)>\oplus <Cl(\mathbf b_j)>$ is not unique. We can have
\begin{displaymath}
\begin{split}
&<Cl(\mathbf b_i)>\oplus <Cl(\mathbf b_j)>=<Cl(\mathbf b_i^\prime)>\oplus <Cl(\mathbf b_j^\prime)>,\\
&\sigma(\mathbf b_i)\simeq \mathbf b_i^\mu,\quad \sigma(\mathbf b_j)\simeq\mathbf  b_j^\mu,\quad
\sigma(\mathbf b_i^\prime)\simeq (\mathbf b_i^\prime)^\mu,\quad \sigma(\mathbf b_j^\prime)
\simeq (\mathbf b_j^\prime)^\mu,\\
& <Cl(\mathbf b_i)>\not\in\{<Cl(\mathbf b_j)>,\quad <Cl(\mathbf b_i^\prime)>,\quad
<Cl(\mathbf b_j^\prime)>\}.\\
\end{split}
\end{displaymath}
\item
Recall that $P_{r_1}(\sigma)\in {\bf F}_p[G]$ is the minimal polynomial 
such that $\mathbf b^{P_{r_1}(\sigma)}\simeq\Z[\zeta]$ with $r_1\leq r_p$.
\item
Recall that $P_{r_1^-}(\sigma)\in {\bf F}_p[G]$ is the minimal polynomial 
such that 
$(\frac{\mathbf b}{ \overline{\mathbf b}})^{P_{r_1}(\sigma)}\simeq\Z[\zeta]$
 with $r_1^-\leq r_p^-$.
\item
We say that the algebraic number  $C\in \Q(\zeta)$ is  singular if $C\Z[\zeta]=\mathbf c^p$ for 
some ideal $\mathbf c$ of $\Q(\zeta)$. We say that $C$ is singular primary 
if $C$ is singular and $C\equiv c^p\modu \pi^p,\quad c\in \Z,\quad c\not\equiv 0\modu p$.
\end{itemize}
%

At fisrt,  a general lemma dealing with congruences on $p$-powers of algebraic numbers of $\Q(\zeta)$.
\begin{lem}$ $\label{l502031}
Let $\alpha,\beta\in\Z[\zeta]$ 
with $\alpha\not\equiv 0 \mod \pi$ 
and $\alpha\equiv\beta\modu \pi$.
Then $\alpha^p\equiv \beta^p  \modu \pi^{p+1}$.
\begin{proof}$ $
Let $\lambda=(\zeta-1)$. 
Then 
$\alpha-\beta\equiv 0 \modu \pi$  implies that 
$\alpha-\zeta^k\beta\equiv 0 \modu \pi$ for $k=0,1,\dots,p-1$. 
Therefore, for all 
$k, \quad 0\leq k \leq p-1$, there exists $a_k\in \N,\quad 0\leq a_k \leq p-1$, such that 
$(\alpha-\zeta^k\beta)\equiv \lambda a_k \modu \pi^2$. 
For another value $l,\quad 0\leq l\leq p-1$, we have, in the same way, 
$(\alpha-\zeta^l\beta)\equiv \lambda a_l \modu \pi^2$,
hence $(\zeta^k-\zeta^l)\beta\equiv \lambda(a_k-a_l) \modu \pi^2$. 
For $k\not=l$ we get $a_k\not= a_l$, because $\pi \| (\zeta^k-\zeta^l)$ and 
because hypothesis $\alpha\not\equiv 0 \modu \pi$ implies that $\beta\not\equiv 0 \modu \pi$. 
Therefore, there exists one and only one $k$ such that 
$(\alpha-\zeta^k\beta)\equiv 0 \modu \pi^2$. Then, we have 
$ \prod_{j=0}^{p-1}(\alpha-\zeta^j\beta)
= (\alpha^p-\beta^p)\equiv 0 \modu \pi^{p+1}$.
\end{proof}
\end{lem}
%
For $i=1,\dots,r_p^-$, to simplify notations in this lemma,  let us note  respectively $\mathbf b, B, C, \mu=u_{2m+1} $ for $\mathbf b_i, B_i, C_i, \mu_i=u_{2m_i+1}$ as defined in the two relations (\ref{e203031}) p. \pageref{e203031} and (\ref{e108191}) p. \pageref{e108191}. 
\begin{lem}\label{l108161}
For $i=1,\dots,r_p^-$, there exists algebraic integers $B\in \Z[\zeta]$  such that
\begin{equation}
\begin{split}
& B\Z[\zeta]=\mathbf b^p,\\
& \sigma(\frac{B}{\overline{B}})
\times(\frac{B}{\overline{B}})^{-\mu}
=(\frac{\alpha}{\overline{\alpha}})^p,\quad\alpha\in\Q(\zeta),
\quad \alpha\Z[\zeta]+\pi=\Z[\zeta],\\
& \sigma(\frac{B}{\overline{B}})\equiv
(\frac{B}{\overline{B}})^{\mu}\modu \pi^{p+1}.
\end{split}
\end{equation}
\begin{proof}$ $
\begin{enumerate}
\item
Observe that we can neglect in this proof the values $\mu=u_{2m}$ such that $\sigma-\mu$ annihilates ideal classes  $\in C_p^+$, because we consider only 
quotients $\frac{B}{\overline{B}}$, with ideal classes  $Cl(\mathbf b)$ in $C_p^-$.
The ideal $\mathbf b^p$ is principal. So let one $\beta\in\Z[\zeta]$ with 
$\beta\Z[\zeta]=\mathbf b^p$.
We have seen in relation (\ref{e107162}) p.\pageref{e107162} that   $\sigma(\mathbf b)\simeq \mathbf b^{\mu}$, therefore there exists 
$\alpha\in \Q(\zeta)$ such that 
$\frac{\sigma(\mathbf b)}{\mathbf b^{\mu}}=\alpha\Z[\zeta]$, also
$\frac{\sigma(\mathbf \beta)}{\mathbf \beta^{\mu}}= \varepsilon\times \alpha^p,\quad \varepsilon\in\Z[\zeta]^*$.
Let $B=\delta^{-1}\times \beta,\quad \delta\in \Z[\zeta]^*$,  for a choice of
the unit $\delta$ that whe shall explicit in the next lines.
We have
\begin{displaymath}
\sigma(\delta\times B)=\alpha^p\times (\delta\times B)^{\mu}\times\varepsilon.
\end{displaymath}
Therefore
\begin{equation}\label{e108162}
\sigma(B)=
\alpha^p\times B^{\mu} \times 
(\sigma(\delta^{-1})\times\delta^{\mu}\times\varepsilon).
\end{equation}
From Kummer's lemma on units, we can write 
\begin{displaymath}
\begin{split}
& \delta=\zeta^{v_1}\times\eta_1,\quad v_1\in \Z,
\quad \eta_1\in\Z[\zeta+\zeta^{-1}]^*,\\ 
& \varepsilon=\zeta^{v_2}\times\eta_2,\quad v_2\in \Z,
\quad \eta_2\in\Z[\zeta+\zeta^{-1}]^*.
\end{split}
\end{displaymath}
Therefore 
\begin{displaymath}
\sigma(\delta^{-1})\times \delta^{\mu}\times\varepsilon=
\zeta^{-v_1 u+v_1\mu+v_2}\times\eta,\quad \eta\in\Z[\zeta+\zeta^{-1}]^*.
\end{displaymath}
From lemma \ref{l103021} p.\pageref{l103021}, we deduce that $\mu\not=u$, therefore there exists one 
$v_1$ with $-v_1 u+v_1\mu+v_2\equiv 0\modu p$.
Therefore, chosing this value $v_1$ for the unit $\delta$, 
\begin{equation}\label{e112212}
\begin{split}
& \sigma(B)=\alpha^p\times B^{\mu}\times\eta,
\quad \alpha\Z[\zeta]+\pi=\Z[\zeta],\quad \eta\in\Z[\zeta+\zeta^{-1}]^*,\\
& \sigma(\overline{B})=\overline{\alpha}^p\times\overline{B}^{\mu}
\times\eta.
\end{split}
\end{equation}
We have $\alpha\equiv \overline{\alpha}\modu \pi$ and we have proved in lemma  \ref{l502031} p.\pageref{l502031} that $\alpha^p\equiv\overline{\alpha}^p\modu \pi^{p+1}$,  which leads to the result.
\end{enumerate}
\end{proof}
\end{lem}
%
\paragraph{$\pi$-adic congruences connected to relative $p$-class group $C_p^-$ :}

For $i=1,\dots,r_p^-$, to simplify notations in this lemma,  let us note  respectively $\mathbf b, B, C, \mu=u_{2m+1} $ for $\mathbf b_i, B_i, C_i, \mu_i=u_{2m_i+1}$ as defined in the two relations (\ref{e203031}) p. \pageref{e203031} and (\ref{e108191}) p. \pageref{e108191}.
\begin{lem}\label{l108171}
For each $i=1,\dots,r_p^-$, there exists singular algebraic integers   $B\in \Z[\zeta]$, such that
\begin{equation}\label{e108191}
\begin{split}
& \mu=u_{2m+1},\quad m\in\N,\quad  1\leq m\leq\frac{p-3}{2},\\
& B\Z[\zeta]=\mathbf b^p,\\
& C=\frac{B}{\overline{B}}\equiv 1
\modu \pi^{2m+1}.
\end{split}
\end{equation}
Then, either $C$ is singular not primary with  $\pi^{2m+1}\ \|\ C-1$  or $C$ is singular primary with $\pi^p\ |\ C-1$.
\begin{proof} $ $
\begin{itemize}
\item  
The definition of $C$ implies that $C\equiv 1\modu \pi$, and so that 
$\sigma(C)\equiv 1\modu\pi$. There exists a natural integer $\nu$ such that $\pi^\nu\ \|\ C-1$, therefore we can write 
\begin{equation}\label{e502031}
\begin{split}
& C\equiv 1+c_0\lambda^\nu\modu \lambda^{\nu+1},\\
& c_0\in\Z,\quad c_0\not\equiv 0\modu p.\\ 
\end{split}
\end{equation}
We have to prove that $\nu<p$ implies that 
$\nu=2m+1$ for the  integer $m<p-1$ verifying $\mu=u_{2m+1}$.
\item
From lemma \ref{l108161} p.\pageref{l108161}, it follows that $\sigma(C)=C^\mu\times\alpha^p$, with 
some $\alpha\in \Q(\zeta)$, and so that 
$1+c_0\sigma(\lambda)^\nu\equiv  (1+\mu c_0\lambda^\nu)\times \alpha^p\modu\pi^{\nu+1}$. 
This congruence implies that $\alpha\equiv 1\modu\pi$ and then, from lemma \ref{l502031},  $\alpha^p\equiv 1\modu\pi^{p+1}$.
Then $1+c_0\sigma(\lambda)^\nu\equiv 1+\mu c_0\lambda^\nu\modu\lambda^{\nu+1}$, and so 
$\sigma(\lambda^{\nu})\equiv \mu\lambda^{\nu}\modu\pi^{\nu+1}$. This implies that $\sigma(\zeta-1)^\nu\equiv \mu\lambda^{\nu}\modu\pi^{\nu+1}$, 
so that $(\zeta^u-1)^\nu\equiv \mu\lambda^{\nu}\modu\pi^{\nu+1}$, so that 
$((\lambda+1)^u-1)^\nu\equiv \mu\lambda^{\nu}\modu\pi^{\nu+1}$ and finally 
$u^\nu\lambda^\nu\equiv \mu\lambda^{\nu}\modu\pi^{\nu+1}$, with simplification $u^\nu-\mu\equiv 0\modu \pi$. Therefore, we have proved that $\nu=2m+1$ or that, 
when $\pi^p\ \not|\ C-1$,  then $\pi^{2m+1}\ \|\ C-1$. 
\end{itemize}
\end{proof}
\end{lem}
%
{\bf Remarks:} 
%
\begin{enumerate}
\item In considering $\mathbf b^{p-1}$ in place of $\mathbf b$, we consider $B^{p-1}$ and $C^{p-1}$ in place of $B$ and $C$, such that we can always assume without loss of generality that $B\equiv C\equiv 1\modu \pi$. We suppose implicitely this normalization in the sequel.
\item
In relation (\ref{e502031}) we can suppose without loss of generality that $c_0=1$ because we can consider $\mathbf b^n$ with $1\leq n \leq p-1$ in place of $\mathbf b$ with $n\times c_0\equiv 1\modu p$. We suppose implicitely this normalization in the sequel.
\end{enumerate}
%
%
As previously, for $i=1,\dots,r_p$, to simplify notations in this lemma,  let us note  respectively $\mathbf b, B, C, \mu=u_{2m+1} $ for $\mathbf b_i, B_i, C_i, \mu_i=u_{2m_i+1}$ as defined in the two relations (\ref{e203031}) p. \pageref{e203031} and (\ref{e108191}) p. \pageref{e108191}. 
In the following lemma, we connect $\pi$-adic congruences on $C-1$ with  
$C=\frac{B}{\overline{B}}$ to some $\pi$-adic congruences on algebraic integer $B$.
%
\begin{thm} ***\label{l203171}
\begin{enumerate}
\item
If the singular number $B$ is not primary, 
there exists  a primary unit $\eta\in\Z[\zeta+\zeta^{-1}]-\{1,-1\}$ and a singular not primary number 
$B^{\prime}= \frac{B^2}{\eta}$, such that 
\begin{equation}\label{e203271}
\begin{split}
& \sigma(B^{\prime})=B^{\prime\mu}\times \alpha^p,\quad \alpha\in\Q(\zeta),\\
& B^{\prime}\Z[\zeta]=\mathbf b^{2p},\quad B^{\prime}\in\Z[\zeta],\\
& \pi^{2m+1}\ \|\ (B^{\prime})^{p-1}-1.\\
\end{split}
\end{equation}
\item
If the singular number B is primary then
\begin{equation}\label{e40721}
\begin{split}
& \sigma(B)=B^{\mu}\times \alpha^p,\quad \alpha\in\Q(\zeta),\\
& B\Z[\zeta]=\mathbf b^{p},\\
& \pi^{p-1}\ |\ B-1.\\
\end{split}
\end{equation}
\end{enumerate}
\begin{proof}
\begin{enumerate}
\item
We have $C\Z[\zeta]=\mathbf b^p$ where the ideal $\mathbf b$ verifies 
$\sigma(\mathbf b)\simeq \mathbf b^{\mu}$ and 
$Cl(\mathbf b)\in C_p^-$. 
From  relation (\ref{e112212}) p. \pageref{e112212} and from 
$Cl(\mathbf b)\in C_p^-$, we can choose $B$ such that 
\begin{displaymath}
\begin{split}
& C=\frac{B}{\overline{B}},\quad B\in\Z[\zeta], \\
& B\overline{B}=\eta\times\gamma^p,\quad \eta\in\Z[\zeta+\zeta^{-1}]^*,
\quad \gamma\in\Q(\zeta),\quad v_\pi(\gamma)=0,\\
& \sigma(B) = B^\mu\times\alpha^p\times \varepsilon,
\quad \mu=u_{2m+1},\quad \alpha\in\Q(\zeta),\quad v_\pi(\alpha)=0,
\quad \varepsilon\in\Z[\zeta+\zeta^{-1}]^*.
\end{split}
\end{displaymath}
We derive  that 
\begin{displaymath}
\begin{split}
& \sigma(B\overline{B})=\sigma(\eta)\times\sigma(\gamma^p)\\
&\sigma(B\overline{B})= (B\overline{B})^\mu\times (\alpha\overline{\alpha})^p\times \varepsilon^2
=\eta^\mu\gamma^{p\mu}\times(\alpha\overline{\alpha})^p\times \varepsilon^2,\\
\end{split}
\end{displaymath}
and so 
\begin{displaymath}
\sigma(\eta)=\eta^\mu\times\varepsilon^2\times\varepsilon_1^p,\quad
\varepsilon_1\in \Z[\zeta+\zeta^{-1}]^*.
\end{displaymath}
We have seen that
\begin{displaymath}
\sigma(B^2)=B^{2\mu}\times  \alpha^{2p}\times \varepsilon^2,
\end{displaymath}
and so 
\begin{displaymath}
\sigma(B^{2})=B^{2\mu}\times\alpha^{2p}\times (\sigma(\eta)\eta^{-\mu}\varepsilon_1^{-p})
\end{displaymath}
which leads to 
\begin{displaymath}
\sigma(\frac{B^2}{\eta})=(\frac{B^{2}}{\eta})^\mu\times \alpha_2^p,\quad \alpha_2^p=\alpha^{2p}\times\varepsilon_1^{-p},\quad 
\alpha_2\in\Q(\zeta),\quad v_\pi(\alpha_2)=0.
\end{displaymath}
Let us note $B^\prime=\frac{B^2}{\eta},\quad B^\prime\in\Z(\zeta),\quad v_\pi(B^\prime)=0$.
We get 
\begin{equation}\label{e203092}
\sigma(B^\prime)=(B^\prime)^\mu\times \alpha_2^p.
\end{equation}
This relation (\ref{e203092}) is similar to hypothesis used to prove  lemma  \ref{l108171} p. \pageref{l108171}. This leads in the same way to 
$B^\prime\equiv d^p\modu \pi^{2m+1},\quad d\in\Z,\quad d\not\equiv 0\modu p$. Therefore $(B^\prime)^{p-1}\equiv 1\modu \pi^{2m+1}$, 
 which achieves  the proof of the first part.
\item
We have 
\begin{equation}\label{e407212}
\begin{split}
& \sigma(B)=B^\mu\times\alpha^p\times\eta,\quad \eta\in\Z[\zeta+\zeta^{-1}]^*\\
& \sigma(\overline{B})=\overline{B}^\mu\times\overline{\alpha}^p\times\eta.\\
\end{split}
\end{equation}
From {\bf simultaneous} application of a Furtwangler theorem, see Ribenboim \cite{rib} (6C) p. 182 and of a Hecke theorem on class field theory, see Ribenboim \cite{rib} (6D) p. 182, it results that  
\begin{equation}\label{e407213}
B\times\overline{B}=\beta^p
\end{equation}
where $\beta\in\Z[\zeta]-\Z[\zeta+\zeta^{-1}]^*$.
From these two relations, it follows that $\eta\in (\Z[\zeta+\zeta^{-1}]^*)^p$, which achieves the proof
of the second part.
\end{enumerate}
\end{proof}
\end{thm}
%
%
\clearpage
\paragraph{On structure of $p$-class group $C_p^-$}
\paragraph{}
In this paragraph, the indexing  of singular primary and of singular not primary $C_i$ with usual previously  defined  meaning of index $i=1,\dots,r_p^+,r_p^++1,\dots,r_p^-$, is used to describe the structure of relative $p$-class group $C_p^-$: we shall show that, with a certain ordering of index $i$, then  $C_i$ are singular primary  for $i=1,\dots,r_p^+$ and 
$C_i$ are singular not primary  for $i=r_p^++1,\dots,r_p^-$.
%
\begin{thm} {*** }\label{t112311}
Let $\mathbf C=C_{1}^{\alpha_{1}}\times \dots\times C_{i}^{\alpha_{i}}\times\dots\times C_{n}^{\alpha_{n}}$ with $\alpha_{i}\in{\bf F}_p^*,
\quad 1\leq n\leq r_p^-$
and with $\mu_{i}=u_{2m_{i}+1}$ pairwise different for $i=1,\dots,n$.
Then $C$ is singular primary if and only if all the $C_{i},\quad i=1,\dots,n$, are all singular primary.
\begin{proof}$ $
\begin{itemize}
\item
If $C_{i},\quad i=1,\dots,n$, are all singular primary, then $C$ is clearly singular primary.
\item
Suppose that $C_{i},\quad i=1,\dots, l$, are not singular primary and that $C_{i},\quad i=l+1,\dots,n$, are singular primary. 
Then, from lemma \ref{l108171} p.\pageref{l108171} and remark following it, 
$\pi^{2m_{i}+1}\|\  C_{i}-1,\quad i=1,\dots,l$, where we suppose, without loss of generality, that
$1<2m_1+1<\dots < 2m_l+1$. Then $\pi^{2m_1+1}\ \|\  C-1$ and so $C$ is not singular primary, contradiction which achieves the proof.
\end{itemize}
\end{proof}
\end{thm}
%
\begin{lem}\label{l202231}
Let $C$ of  relation (\ref{e108191}) p. \pageref{e108191}. If $C$ is not singular primary, then 
\begin{displaymath}
C\equiv 1+ V(\mu)\modu \pi^{p-1},\quad \mu=u_{2m+1}
\quad V(\mu)\in\Z[\zeta],
\end{displaymath}
where $V(\mu)\modu p$ depends only on $\mu$ with $\pi^{2m+1}\|\  V(\mu)$.
\begin{proof}
The congruence $\sigma(C)\equiv C^\mu\modu\pi^{p+1}$ and the normalization $C\equiv 1+\lambda^{2m+1}\modu\pi^{2m+2}$ explained in remark following
lemma \ref{l108171} p. \pageref{l108171} implies the result.
\end{proof}
\end{lem}
%
\begin{thm}{ *** }\label{l203031}
Let $C_{1}, C_{2}$ singular not primary defined with relation (\ref{e108191}) p. \pageref{e108191}.
If $\mu_{1}=\mu_{2}$ then $C_{1}\times C_{2}^{-1}$ is singular  primary.
\begin{proof}
Let $\mu_{1}=\mu_{2}=\mu=u_{2m+1}$. Therefore $\sigma(C_{1})\equiv C_{1}^\mu\modu\pi^{p+1}$ and 
$\sigma(C_{2})\equiv C_{2}^\mu\modu\pi^{p+1}$.
From previous lemma \ref{l202231} p. \pageref{l202231} we get 
\begin{displaymath}
\begin{split}
& C_{1}=1+V(\mu)+p W_1,\quad W_1\in\Q(\zeta), \quad v_\pi(V(\mu))\geq 2m+1,
\quad v_\pi(W_1)\geq 0,\\
& C_{2}=1+V(\mu)+p W_2,\quad W_2\in\Q(\zeta), \quad v_\pi(V(\mu))\geq 2m+1,
\quad v_\pi(W_2)\geq 0,\\
& \pi^{2m+1}\ \|\ V(\mu),\\
\end{split}
\end{displaymath}
Elsewhere, $C_{1}, C_{2}$ verify
\begin{displaymath}
\begin{split}
& \sigma(C_{1})\equiv C_{1}^\mu\modu\pi^{p+1},\\
& \sigma(C_{2})\equiv C_{2}^\mu\modu\pi^{p+1},
\end{split}
\end{displaymath}
which leads to 
\begin{displaymath}
\begin{split}
& 1+\sigma(V(\mu))+p \sigma(W_1)\equiv 
1+ A(\mu)+ p \mu W_1 \modu \pi^{p+1},\\
& 1+\sigma(V(\mu))+p \sigma(W_2)\equiv 
1+ A(\mu)+ p \mu  W_2 \modu \pi^{p+1},
\end{split}
\end{displaymath}
where $A(\mu)\in \Q(\zeta),\quad v_\pi(A(\mu)\geq 0$ depends only on $\mu$.
By difference, we get 
\begin{displaymath}
p(\sigma(W_1-W_2))\equiv p  \mu (W_1- W_2)\modu \pi^{p+1},
\end{displaymath}
which implies that
\begin{displaymath}
\sigma(W_1-W_2)\equiv   \mu (W_1- W_2)\modu \pi^{2}.
\end{displaymath}
Let $W_1-W_2= a\lambda+b,\quad a,b\in\Z,\quad \lambda=\zeta-1$.
The previous relation implies that $b(1-\mu)\equiv 0\modu p$ and so 
that $a\sigma(\lambda)+b\equiv \mu a\lambda+\mu b\modu\pi^2$, 
and so that $b\equiv 0\modu p$, because $\mu\not=1$.
Thus $W_1-W_2\equiv 0\modu \pi$
and finally $C_1\equiv C_2\modu \pi^p$ and also
$C_1 C_2^{-1}\equiv 1\modu \pi^p$ and $C_1 C_2^{-1}$ is singular primary.
\end{proof}
\end{thm}
%
\begin{cor}\label{c301031}
Let $C_{1},\dots, C_{\nu},\quad 1\leq \nu\leq r_p^-$, singular not primary, defined by relation (\ref{e108191}) p. \pageref{e108191}.
\begin{enumerate}
\item
If $\mu_{1}=\dots=\mu_{\nu}=\mu$ then $C_{1}^\prime=C_{1}\times C_{\nu}^{-1},\dots, 
C_{\nu-1}^\prime=C_{\nu-1}\times  C_{\nu}^{-1}$ are singular primary.
\item
In term of ideals,  it implies that
\begin{displaymath}
\oplus_{i=1}^\nu <Cl(\mathbf b_{i})>=\oplus_{i=1}^{\nu-1} <Cl(\mathbf b_{i}\mathbf b_{\nu}^{-1})>
\oplus <Cl(\mathbf b_{\nu})>,
\end{displaymath}
where $\sigma(\mathbf b_i\mathbf b_\nu^{-1})\simeq (\mathbf b_i\mathbf b_\nu^{-1})^{\mu}$
\end{enumerate}
\begin{proof}$ $
\begin{enumerate}
\item
Immediate consequence of theorem \ref{l203031} p. \pageref{l203031}.
\item
$\oplus_{i=1}^\nu <Cl(\mathbf b_{i})>$ is a $p$-group of rank $\nu$.
$\oplus_{i=1}^{\nu-1} <Cl(\mathbf b_{i}\mathbf b_{\nu}^{-1})>$ is a $p$-group of rank $\nu-1$.
$<Cl(\mathbf b_\nu)>$ is a $p$-group of rank $1$.
\end{enumerate}
\end{proof}
\end{cor}
%
\paragraph{Remark:} $ $
It follows that, when $\mu_{1}=\dots=\mu_{\nu}=\mu$, we can suppose without loss of generality, with usual meaning of indexing $i=1,\dots,r_p^-$, that  the representants
$C_1,\dots, C_{\nu-1}$  chosen are singular primary.
%
%
\begin{thm}{ *** On structure of $p$-class group $C_p^-$. }\label{t203021}

Let $\mathbf b_i$ be the ideals defined in relation (\ref{e203031}) p. \pageref{e203031}.
Let $C_p^{-}=\oplus _{i=1}^{r_p^-} <Cl(\mathbf b_i)>$.
Let $C_i=\frac{B_i}{\overline{B}_i},\quad B_i\Z[\zeta]=\mathbf b_i^p,
\quad Cl(\mathbf b_i)\in C_p^-,
\quad i=1,\dots,r_p^-$,  where
$B_i$ is defined in relation (\ref{e112212}) p. \pageref{e112212}. 
Let $r_1^-$ be the degree of the minimal polynomial $P_{r_1^-}(\sigma)$ defined in relation (\ref{e203031}) p. \pageref{e203031}.
With a certain ordering of $C_i,\quad i=1,\dots,r_p^-$,
\begin{enumerate}
\item
$C_i$ are  singular primary for $i=1,\dots, r_p^+$, and $C_i$ are  singular not primary for $i=r_p^+,\dots, r_p^-$.
\item
\begin{enumerate}
\item
If $j>i\geq r_p^++1$  then $\mu_j\not=\mu_i$.
\item
If $\mu_i=\mu_j$ then $j<i\leq r_p^+$.
\end{enumerate}
\item
$r_p^--r_p^+\leq r_1^-\leq r_p^-$.
\end{enumerate}
\begin{proof}$ $
\begin{enumerate}
\item
It is an application  of a theorem of Furwangler, see Ribenboim \cite{rib} (6C) p. 182 and of a theorem of Hecke, see Ribenboim \cite{rib} (6D) p. 182.
\item 
See lemma \ref{l203031} p. \pageref{l203031}
\item
Apply corollary \ref{c301031} p. \pageref{c301031}.
\end{enumerate}
\end{proof}
\end{thm}
%
%
\subsubsection{$\pi$-adic congruences  connected  to  $p$-class group $C_p^+$}
%
For $i=r_p^++1,\dots,r_p$, to simplify notations in this lemma,  let us note  respectively $\mathbf b, B$ for ideal $\mathbf b_i$ and algebraic integer $B_i$, as defined in the two relations (\ref{e203031}) p. \pageref{e203031} and (\ref{e108191}) p. \pageref{e108191}.
%
\begin{thm} *** \label{l203272}

Let the ideals $\mathbf b$, such that 
$Cl(\mathbf b)\in C_p^+$ defined in relation (\ref{e203031}) p. \pageref{e203031}. There exists 
$B\in\Z[\zeta]$ such that:
\begin{equation}\label{e205041}
\begin{split}
& \mu=u_{2n}, \quad 1\leq n\leq \frac{p-3}{2}\\
& \sigma(\mathbf b)\simeq \mathbf b^{\mu},\\
& B\Z[\zeta]=\mathbf b^p,\\
& \sigma(B)=B^{\mu}\times\alpha^p,\quad\alpha\in\Q(\zeta),\\
& B\equiv 1\modu \pi^{2n}.\\
\end{split}
\end{equation}
\begin{proof}
Similarly to relation (\ref{e112212}) p. \pageref{e112212}, there exists $B$ with $B\Z[\zeta]=\mathbf b^p$ such that 
\begin{displaymath}
\sigma(B)=B^{\mu}\times\alpha^p\times \eta,
\quad \alpha\in\Q(\zeta),\quad \eta\in\Z[\zeta+\zeta^{-1}]^*.
\end{displaymath}
From relation (\ref{e201274}) p. \pageref{e201274}, independant forward reference in section dealing of 
unit group $\Z[\zeta+\zeta^{-1}]^*$, we can write
\begin{displaymath}
\begin{split}
& \eta= \eta_1^{\lambda_1}\times (\prod_{j=2}^{N} \eta_{j}^{\lambda^j}),
\quad \lambda_j\in {\bf F}_p,\quad 1\leq  N<\frac{p-3}{2},\\
& \sigma(\eta_1)=\eta_1^{\mu}\times\beta_1^p,
\quad \eta_1,\beta_1\in\Z[\zeta+\zeta^{-1}]^*,\\
& \sigma(\eta_{j})=\eta_{j}^{\nu_{j}}\times\beta_j^p,
\quad \eta_j,\beta_j\in\Z[\zeta+\zeta^{-1}]^*, \quad j=2,\dots,N,\\
& 2\leq j<j^\prime\leq N\Rightarrow \nu_j\not=\nu_{j^\prime},\\
& \nu_j\not=\mu,\quad j=2,\dots,N.\\
\end{split}
\end{displaymath}
Let us note 
\begin{displaymath}
E=\eta_1^{\lambda_1},\quad 
U=\prod_{j=2}^{N} \eta_j^{\lambda_j}.
\end{displaymath}
Show that there exists $V\in\Z[\zeta+\zeta^{-1}]^*$ such that 
\begin{displaymath}
\sigma(V)\times V^{-\mu}=U^{-1}\times \varepsilon^p,\quad \varepsilon\in\Z[\zeta+\zeta^{-1}]^*
\end{displaymath}
Let us suppose that $V$ is of form  $V=\prod_{j=2}^{N} \eta_j^{\rho_j}$.
Then, it suffices  that
\begin{displaymath}
\eta_j^{\rho_j \nu_j}\times \eta_j^{-\rho_j\mu}
=\eta_j^{-\lambda_j}\times \varepsilon_j^p,
\quad \varepsilon_j\in\Z[\zeta+\zeta^{-1}]^*,\quad j=2,\dots,N.
\end{displaymath}
It suffices that
\begin{displaymath}
\rho_j\equiv \frac{-\lambda_j}{\nu_j-\mu}\modu p,\quad j=2,\dots,N,
\end{displaymath}
which is possible, because $\nu_j\not\equiv \mu,\quad j=2,\dots, N$.
Therefore, for $B^\prime=B\times V$, we get  $B=B^\prime V^{-1}$ and so
\begin{displaymath}
\sigma(B)=\sigma(B^\prime V^{-1})=B^\mu\alpha^p\eta=(B^\prime V^{-1})^\mu\alpha^p\eta=(B^\prime V^{-1})^{\mu}\times\alpha^p\times E\times U,
\end{displaymath}
so 
\begin{displaymath}
\sigma(B^\prime )=(B^\prime)^{\mu}\sigma(V) V^{-\mu}\times\alpha^p\times E\times U,
\end{displaymath}
so 
\begin{displaymath}
\sigma(B^\prime )=(B^\prime)^{\mu}(U^{-1}\varepsilon^p\times\alpha^p\times E\times U),
\end{displaymath}
so we get simultaneously 
\begin{displaymath}
\begin{split}
& \sigma(B^\prime)=(B^\prime)^{\mu}\times \alpha^p\times\varepsilon^p\times E,\quad \alpha\in\Q(\zeta),\\
& \sigma(E)=E^{\mu}\times \varepsilon_1^p,\quad \varepsilon_1\in\Z[\zeta+\zeta^{-1}]^*.\\
\end{split}
\end{displaymath}
Show that 
\begin{equation}\label{e407214}
\sigma(B^\prime)=B^{\prime\mu}\times \alpha^{\prime p}.
\end{equation}
\begin{enumerate}
\item
If $E\in (\Z[\zeta+\zeta^{-1}]^*)^p$, then we get
\begin{equation}\label{e203281}
\sigma(B^\prime)=B^{\prime\mu}\times \alpha^{\prime p}.
\end{equation}
\item
If $E\not\in(\Z[\zeta+\zeta^{-1}]^*)^p$ then  by conjugation $\sigma$,
\begin{displaymath}
\begin{split}
& \sigma(B^\prime)=(B^\prime)^{\mu}\times \alpha_1^{ p}\times E,\quad \alpha_1\in\Q(\zeta),\\
& \sigma^2(B^\prime)=\sigma(B^\prime)^{\mu}\times  E^{\mu}\times b^p,
\quad b\in\Q(\zeta),\\
\end{split}
\end{displaymath}
so gathering these relations
\begin{displaymath}
\begin{split}
& \sigma(B^\prime)^{\mu}=(B^\prime)^{\mu^2}\times \alpha_1^{p\mu}\times E^{\mu},\\
& \sigma^2(B^\prime)=\sigma(B^\prime)^{\mu}\times  E^{\mu}\times b^p,\\
\end{split}
\end{displaymath}
and so
\begin{displaymath}
c^p\times \sigma(B^\prime)^{\mu}(B^\prime)
^{-\mu^2}=\sigma^2(B^\prime)\sigma(B^\prime)^{-\mu},
\quad c\in\Q(\zeta)
\end{displaymath}
which leads to
\begin{displaymath}
c^pB^{\prime\mu\sigma}(B^\prime)^{-\mu^2}=B^{\prime\sigma^2}(B^\prime)^{-\mu\sigma},
\end{displaymath}
and so
\begin{displaymath}
(B^\prime)^{(\sigma-\mu)^2}=c^p,\quad c\in\Q(\zeta).
\end{displaymath}
Elsewhere $(B^\prime)^{\sigma^{p-1}-1}=1$, so
\begin{displaymath} 
(B^\prime)^{gcd((\sigma^{p-1}-1, 
(\sigma-\mu)^2})= (B^\prime)^{\sigma-\mu}=\alpha_3^{ p},
\quad \alpha_3\in\Q(\zeta),
\end{displaymath}
and so $\sigma(B^\prime)=(B^\prime)^{\mu}\times\alpha_3^{ p}$.
\end{enumerate}
The end of proof is similar to 
proof of previous lemma \ref{l203171} p. \pageref{l203171}.
\end{proof}
\end{thm}
%
%
\subsubsection{$\pi$-adic congruences  connected  to  $p$-class group $C_p$}
%
Let us consider the  ideals $\mathbf b_i,\quad i=1,\dots,r_p,$ defined in relation 
(\ref{e203031}) p. \pageref{e203031}. Then $Cl(\mathbf b_i)\in C_p$. 
From  theorem \ref{l203171} p. \pageref{l203171}, and theorem \ref{l203272} p. \pageref{l203272},
we can  choose the corresponding singular primary number $B_i$ with
$B_i\Z[\zeta]=\mathbf b_i^p$; then $\sigma(B_i)=B_i^\mu\times\alpha_i^p,\quad \alpha_i\in\Q(\zeta),\quad
\mu_i=u_{m_i},\quad 1\leq m_i\leq p-2$. Observe that if $m_i$ is odd then $Cl(\mathbf b_i)\in C_p^-$ and if 
$m$ is even then $Cl(\mathbf b_i)\in C_p^+$.

%

%
The next important theorem {\bf summarize} for all the $p$-class group $C_p$ the  previous theorems 
\ref{l108171} p. \pageref{l108171} and \ref{l203171} p. \pageref{l203171} for relative $p$-class group $C_p^-$ and \ref{l203272} p. \pageref{l203272} for $p$-class group $C_p^+$ 
and give explicit $\pi$-adic congruences connected to $p$-class group of $\Q(\zeta)$.
%
\begin{thm} { *** }{ $\pi$-adic structure of $p$-class group $C_p$}\label{t203281}

Let the ideals $\mathbf b_i,\quad i=1,\dots,r_p$, such that 
$Cl(\mathbf b_i)\in C_p$ and defined by relation (\ref{e203031}) p. \pageref{e203031}. 
Then, there exists singular algebraic integers $B_i\in\Z[\zeta],\quad i=1,\dots,r_p$, such that
\begin{equation}\label{e205042}
\begin{split}
& \mu_i=u_{m_i},\quad 1\leq m_i\leq p-2,\quad i=1,\dots,r_p,\\
& \sigma(\mathbf b_i)\simeq \mathbf b_i^{\mu_i},\\
& B_i\Z[\zeta]=\mathbf b_i^p,\\
& \sigma(B_i)=B_i^{\mu_i}\times\alpha_i^p,\quad \alpha_i\in\Q(\zeta),\\
& B_i\equiv 1\modu \pi^{m_i}.\\
\end{split}
\end{equation}
Moreover, with a certain reindexing of  $i=1,\dots,r_p$:
\begin{enumerate}
\item
The $r_p^+$ {\bf singular}  integers  $B_i,\quad i=1,\dots,r_p^+$, corresponding to $\mathbf b_i\in C_p^-$ 
 are  {\bf primary} with $\pi^p\ |\ B_i-1$. 
\item
The $r_p^--r_p^+$ {\bf singular}  integers  $B_i,\quad i=r_p^++1,\dots,r_p^-$, corresponding to $\mathbf b_i\in C_p^-$ 
are  {\bf not primary}
and verify $\pi^{m_i}\ \|\  (B_i-1)$.
\item
The $r_p^+$ {\bf singular}  numbers  $B_i,\quad i=r_p^-+1,\dots,r_p$, corresponding 
to $\mathbf b_i\in C_p^+$  are  {\bf primary} or {\bf not primary} (without being able to say more)
and verify $\pi^{m_i}\ |\  (B_i-1)$.
\end{enumerate} 
\begin{proof}$ $
\begin{enumerate}
\item
For the case $C_p^-$, apply lemmas \ref{l203171} p. \pageref{l203171} and theorem \ref{t203021} p. \pageref{t203021}.
Toward this result, observe also that if $B_i$ is not primary, then $\pi^{2m_i+1}\ \|\ (B_i^\prime)^{p-1}-1$ and so 
$\pi^{p-1}\ \not|\ (B_i^\prime)^{p-1}-1$  and $C^\prime=(\frac{B^\prime}{\overline{B^\prime}})^{p-1}$ is not singular primary, therefore 
$B_i^\prime$ primary $\Leftrightarrow  C_i^\prime=\frac{B_i^\prime}{\overline{B}_i^\prime}$ primary. 
\item
For the case $C_p^+$ apply the theorem \ref{l203272}
\end{enumerate}
\end{proof}
\end{thm}
%

%
\paragraph {The case $\mu=u_{2m+1}$ with $2m+1>\frac{p-1}{2}$}\label{s201111}
\paragraph{ }
In the next lemma \ref{l201091} p. \pageref{l201091} and theorem \ref{l202211} p. \pageref{l202211},  we shall investigate more deeply the consequences of the congruence 
$C\equiv 1\modu \pi^{2m+1}$ of lemma \ref{l108171} p. \pageref{l108171} when $2m+1> \frac{p-1}{2}$.
%
\begin{lem}\label{l201091}
Let $C$ with $\mu=u_{2m+1},\quad 2m+1>\frac{p-1}{2}$ written in the form:
\begin{displaymath}
\begin{split}
& C=1+\gamma+\gamma_0\zeta+\gamma_1 \zeta^u
+\dots + \gamma_{p-3}\zeta^{u_{p-3}},\\
& \gamma\in\Q,\quad v_p(\gamma)\geq 0,\quad  \gamma_i\in\Q,\quad  
v_p(\gamma_i)\geq 0,\quad i=0, \dots, p-3,\\
& \gamma+\gamma_0\zeta+\gamma_1\zeta^u+\dots + \gamma_{p-3}\zeta^{u_{p-3}}\equiv 0\modu \pi^{2m+1},
\quad 2m+1> \frac{p-1}{2}.
\end{split}
\end{displaymath}
Then $C$ verifies  the congruences 
\begin{displaymath}
\begin{split}
& \gamma\equiv -\frac{\gamma_{p-3}}{\mu-1}\modu p,\\
& \gamma_0\equiv -\mu^{-1} \times \gamma_{p-3}\modu p,\\
& \gamma_1\equiv -(\mu^{-2}+\mu^{-1})\times \gamma_{p-3}\modu p,\\
&\vdots\\
& \gamma_{p-4}\equiv -(\mu^{-(p-3)}+\dots+\mu^{-1})\times \gamma_{p-3}\modu p.
\end{split}
\end{displaymath}
\begin{proof}
We have seen in lemma \ref{l108161} p. \pageref{l108161} that $\sigma(C)\equiv C^{\mu}\modu\pi^{p+1}$.
From $2m+1 > \frac{p-1}{2}$ we derive that
\begin{displaymath}
C^\mu\equiv 1+\mu\times(\gamma+\gamma_0\zeta+\gamma_1\zeta^u
+\dots +\gamma_{p-3} \zeta^{u_{p-3}})\modu\pi^{p-1}.
\end{displaymath}
Elsewhere, we get by conjugation
\begin{equation}\label{e201111}
\sigma(C)=1+\gamma+\gamma_0\zeta^u+\gamma_1\zeta^{u_2}
+\dots +\gamma_{p-3} \zeta^{u_{p-2}}.
\end{equation} 
We have the identity
\begin{displaymath}
\gamma_{p-3} \zeta^{u_{p-2}}=-\gamma_{p-3}-\gamma_{p-3}\zeta
-\dots -\gamma_{p-3}\zeta^{u_{p-3}}.
\end{displaymath}
This leads to 
\begin{displaymath}
\sigma(C)=1+\gamma-\gamma_{p-3}-\gamma_{p-3} 
\zeta+(\gamma_0-\gamma_{p-3})\zeta^u
+\dots+(\gamma_{p-4}-\gamma_{p-3})\zeta^{u_{p-3}}.
\end{displaymath}
Therefore, from the congruence $\sigma(C)\equiv C^{\mu}\modu \pi^{p+1}$ we get
the congruences in the basis $1,\zeta,\zeta^{u},\dots,\zeta^{u_{p-3}}$,
\begin{displaymath}
\begin{split}
& 1+\mu \gamma\equiv 1+\gamma-\gamma_{p-3}\modu p,\\
& \mu \gamma_0\equiv -\gamma_{p-3}\modu p,\\
& \mu \gamma_1\equiv \gamma_0-\gamma_{p-3}\modu p,\\
& \mu \gamma_2\equiv \gamma_1-\gamma_{p-3}\modu p,\\
& \vdots\\
& \mu \gamma_{p-4}\equiv \gamma_{p-5}-\gamma_{p-3}\modu p,\\
& \mu \gamma_{p-3}\equiv \gamma_{p-4}-\gamma_{p-3}\modu p.
\end{split}
\end{displaymath}
From these congruences, we get $\gamma\equiv -\frac{\gamma_{p-3}}{\mu-1}\modu p$ and 
$\gamma_0\equiv -\mu^{-1} \gamma_{p-3}\modu p$ and then
$\gamma_1\equiv \mu^{-1}(\gamma_0-\gamma_{p-3})\equiv \mu^{-1}(-\mu^{-1} 
\gamma_{p-3}-\gamma_{p-3})
\equiv -(\mu^{-2}+\mu^{-1})\gamma_{p-3}\modu p$ 
and $\gamma_2
\equiv\mu^{-1}(\gamma_1-\gamma_{p-3})
\equiv \mu^{-1}(-(\mu^{-2}+\mu^{-1})\gamma_{p-3}-\gamma_{p-3})
\equiv -(\mu^{-3}+\mu^{-2}+\mu^{-1})\gamma_{p-3}\modu p$ and so on.
\end{proof}
\end{lem}
%
%
%
The next theorem gives an explicit  important formulation  of $C$ when $2m+1>\frac{p-1}{2}$.  
\begin{thm}{ *** }\label{l202211}
Let $\mu=u_{2m+1},\quad p-2\geq 2m+1>\frac{p-1}{2}$, corresponding to $C$ defined in lemma \ref{l108171} p. \pageref{l108171}, so $\sigma(C)\equiv C^\mu\modu\pi^{p+1}$.
Then $C$ verifies the formula: 
\begin{equation}\label{e202211}
C\equiv 1-\frac{\gamma_{p-3}}{\mu-1}\times
(\zeta+\mu^{-1}\zeta^u+\dots+\mu^{-(p-2)}\zeta^{u_{p-2}}) \modu \pi^{p-1}.
\end{equation}
\begin{proof}
From definition of $C$, setting $C=1+V$, we get :
\begin{displaymath}
\begin{split}
& C=1+V,\\
& V=\gamma+\gamma_0\zeta+\gamma_1\zeta^u+\dots+\gamma_{p-3}\zeta^{u_{p-3}},\\
& \sigma(V)\equiv \mu\times V\modu \pi^{p+1}.
\end{split}
\end{displaymath}
Then, from lemma \ref{l201091} p. \pageref{l201091},  we obtain the relations 
\begin{displaymath}
\begin{split}
& \mu=u_{2m+1},\\
& \gamma \equiv -\frac{\gamma_{p-3}}{\mu-1}\modu p,\\
& \gamma_0\equiv -\mu^{-1} \times \gamma_{p-3}\modu  p,\\
& \gamma_1\equiv -(\mu^{-2}+\mu^{-1})\times \gamma_{p-3}\modu p,\\
&\vdots\\
& \gamma_{p-4}\equiv -(\mu^{-(p-3)}+\dots+\mu^{-1})\times \gamma_{p-3}\modu  p,\\
& \gamma_{p-3}\equiv -(\mu^{-(p-2)}+\dots+\mu^{-1})\times \gamma_{p-3}\modu  p.
\end{split}
\end{displaymath}
From these relations we get 
\begin{displaymath}
V\equiv  -\gamma_{p-3}\times
(\frac{1}{\mu-1}+\mu^{-1}\zeta+(\mu^{-2}+\mu^{-1})\zeta^{u}+\dots+
(\mu^{-(p-2)}+\dots+\mu^{-1})\zeta^{u_{p-3}})\modu p.
\end{displaymath}
Then 
\begin{displaymath}
V\equiv-\gamma_{p-3}\times (\frac{1}{\mu-1}+
\mu^{-1}(\zeta+(\mu^{-1}+1)\zeta^u+\dots+(\mu^{-(p-3)}+\dots+1)
\zeta^{u_{p-3}}))\modu p.
\end{displaymath}
Then
\begin{displaymath}
V\equiv -\gamma_{p-3}\times (\frac{1}{\mu-1}+
\mu^{-1}(
\frac{(\mu^{-1}-1)\zeta+(\mu^{-2}-1)\zeta^u+\dots
+(\mu^{-(p-2)}-1)\zeta^{u_{p-3}}}{\mu^{-1}-1}))
\modu p.
\end{displaymath}
Then
\begin{displaymath}
V\equiv -\gamma_{p-3}\times (\frac{1}{\mu-1}+
\mu^{-1}(
\frac{\mu^{-1}\zeta+\mu^{-2}\zeta^u+\dots
+\mu^{-(p-2)}\zeta^{u_{p-3}}
-\zeta-\zeta^u-\dots -\zeta^{u_{p-3}}}{\mu^{-1}-1}))
\modu p.
\end{displaymath}
Then  $-\zeta-\zeta^u-\dots-\zeta^{u_{p-3}}=1+\zeta^{u_{p-2}}$ and $\mu^{-(p-1)}\equiv 1\modu p$
implies that 
\begin{displaymath}
V\equiv -\gamma_{p-3}\times (\frac{1}{\mu-1}+
\mu^{-1}(
\frac{1+\mu^{-1}\zeta+\mu^{-2}\zeta^u+\dots
+\mu^{-(p-2)}\zeta^{u_{p-3}}+\mu^{-(p-1)}\zeta^{u_{p-2}}}{\mu^{-1}-1}))
\modu p.
\end{displaymath}
Then
\begin{displaymath}
V\equiv -\gamma_{p-3}\times (\frac{1}{\mu-1})\times 
(1- (1+\mu^{-1}\zeta+\mu^{-2}\zeta^u+\dots
+\mu^{-(p-2)}\zeta^{u_{p-3}}+\mu^{-(p-1)}\zeta^{u_{p-2}}))
\modu p.
\end{displaymath}
Then $\frac{1}{\mu-1}+\frac{\mu^{-1}}{\mu^{-1}-1}=0$ and so 
\begin{displaymath}
V\equiv -\gamma_{p-3}\times (\frac{\mu^{-1}}{\mu-1})\times 
(\zeta+\mu^{-1}\zeta^u+\dots
+\mu^{-(p-3)}\zeta^{u_{p-3}}+\mu^{-(p-2)}\zeta^{u_{p-2}})
\modu p.
\end{displaymath}
\end{proof}
\end{thm}
%
%
%
\clearpage
\section{On  structure of the $p$-unit group of the cyclotomic field $\Q(\zeta)$}\label{s205231}\label{s210163}
Let us consider the results obtained in subsection \ref{s108311} p.\pageref{s108311} for the action of $Gal(\Q(\zeta)/\Q)$ on $C_p^-$. In the present section, we assert that this approach can be partially translated {\it mutatis mutandis} to the study of the  $p$-group of units of $\Q(\zeta)$
\begin{displaymath} 
F=\{\Z[\zeta+\zeta^{-1}]^*/(\Z[\zeta+\zeta^{-1}]^*)^p\}/<-1>. 
\end{displaymath}
%
This section contains:
\begin{itemize}
\item
Some general definitions and properties of the $p$-unit group $F$.
\item
Some $\pi$-adic congruences strongly connected to structure of the $p$-unit-group $F$. These congruences are of the same kind that those found in previous chapter for $p$-class group $C_p$.
\end{itemize}
\subsection{Definitions and preliminary results}
\begin{itemize}
\item
When  $h^-\equiv 0\modu p$,  from Hilbert class field theory, there exists  {\it primary} units 
$\eta\in\Z[\zeta+\zeta^{-1}]^*$, so such that
\begin{equation}\label{e203171}
\begin{split}
& \eta\equiv d^p\modu p,\quad d\in\Z,\quad d\not\equiv 0,\\ 
& \sigma(\eta)=\eta^\mu\times\varepsilon^p,
\quad \varepsilon\in\Z[\zeta+\eta^{-1}]^*.
\end{split}
\end{equation} 
\item
The group $\Z[\zeta+\zeta^{-1}]^*$ is a free group of rank $\frac{p-1}{2}$. It contains the subgroup 
$\{-1,1\}$ of rank $1$.
For all $\eta\in\Z[\zeta+\zeta^{-1}]^*- \{-1,1\}$
\begin{displaymath}
\eta\times\sigma(\eta)\times\dots\times\sigma^{(p-3)/2}(\eta)=\pm 1.
\end{displaymath}
Therefore, for each unit $\eta\in\Z[\zeta+\zeta^{-1}]^*$,  there exists a {\bf minimal} 
$r_\eta\in\N,\quad r_\eta\leq\frac{p-3}{2}$, such that
\begin{equation}\label{e203051}
\begin{split}
& \eta\times\sigma(\eta)^{l_1}\times\dots\times\sigma^{r_\eta}(\eta)^{l_{r_\eta}}
=\varepsilon^p,\quad \varepsilon\in\Z[\zeta+\zeta^{-1}]^*,\\
& 0\leq l_i\leq p-1,\quad i=1,\dots,r_\eta,\quad l_{r_\eta}\not= 0.
\end{split}
\end{equation}
\item
Let us define an equivalence on units of $\Z[\zeta+\zeta^{-1}]^*$ :
$\eta, \eta^\prime\in\Z[\zeta+\zeta^{-1}]^*$ are said equivalent if there exists  
$\varepsilon\in\Z[\zeta+\zeta^{-1}]^*$ such that 
$\eta^\prime= \eta\times \varepsilon^p$. Let us denote $E(\eta)$ the equivalence class of $\eta$. 
\item
We have $E(\eta_a\times\eta_b)=E(\eta_a)\times E(\eta_b)$; the set of class $E(\eta)$ is a group.  The group $<E(\eta)>$ generated by $E(\eta)$ is cyclic of order $p$.
\item 
Observe that this equivalence is consistent with conjugation
$E(\sigma(\eta))=\sigma(E(\eta))$.
\item
The group $F=\{\Z[\zeta+\zeta^{-1}]^*/(\Z[\zeta+\zeta^{-1}]^*)^p\}/<-1>$  so defined is a group of 
rank $\frac{p-3}{2}$, see for instance Ribenboim \cite{rib} p 184 line 14.
\item
Similarly to relation (\ref{e107162}) p.\pageref{e107162}, there exists 
$\eta\in \Z[\zeta+\zeta^{-1}]^*$ such that 
\begin{equation}\label{e203041}
\begin{split}
& E(\eta) = E(\eta_1)\times\dots \times E(\eta_{(p-3)/2}),\\
&E(\sigma(\eta_i))= E(\eta_i^{\mu_i}),\quad i=1,\dots,\frac{p-3}{2},\\
& \mu_i\in \N,\quad 1<\mu_i\leq p-1,\\
& F= <E(\eta_1)>\oplus\dots\oplus <E(\eta_{(p-3)/2})>,
\end{split}
\end{equation}
where $F$ is seen as a ${\bf F}_p[G]$-module of dimension $\frac{p-3}{2}$.
\item
For each unit $\eta$, there is a minimal polynomial $P_{r_\eta}(V)=\prod_{i=1}^{r_\eta} (V-\mu_i)$
where $r_\eta\leq \frac{p-3}{2}$,
such that
\begin{equation}\label{e203042}
\begin{split} 
& E(\eta)^{P_{r_\eta}(\sigma)}=E(1),\\
& 1\leq i < j \leq r_\eta\Rightarrow \mu_i\not=\mu_j.
\end{split}
\end{equation}
\item
Let $\beta\in\Z[\zeta+\zeta^{-1}]^*$.
Observe that if $E(\beta)=E(\sigma(\beta))$ then $E(\sigma^2(\beta))=E(\beta)$ and so $E(1)=E(\beta^{p-1})$ and $E(\beta)=1$.
\item 
Recall that a unit $\beta\in\Z[\zeta+\zeta^{-1}]^*$ is said primary 
if $\beta\equiv b^p\modu\pi^{p+1},\quad b\in\Z$.
\end{itemize}
%
\begin{lem}\label{l106033}
Let $\beta\in\Z[\zeta+\zeta^{-1}]^*-(\Z[\zeta+\zeta^{-1}]^*)^p$.
Then the minimal polynomial $P_{r_{\beta}}(V)$ is of the form
\begin{displaymath}
P_{r_{\beta}}(V)=\prod_{i=1}^{r_{\beta}} (V-u_{2 m_i}),
\quad 1\leq m_i\leq \frac{p-3}{2},\quad r_\beta>0.
\end{displaymath}
\begin{proof}
There exists  $\eta_1\in\Z[\zeta+\zeta^{-1}]^*,\quad E(\eta_1)\not=E(1)$,  with  $E(\eta_1)^{\sigma-\mu_1}=E(1)$.
Suppose that $\mu_1^{(p-1)/2}=-1$ and search for a contradiction:
we have $E(\eta_1)^{\sigma-\mu_1}=E(1)$, therefore 
$E(\eta_1)^{\sigma^{(p-1)/2}-\mu_1^{(p-1)/2}}=E(1)$; but, from
$\eta_1\in\Z[\zeta+\zeta^{-1}]^*$, we get $\eta_1^{\sigma^{(p-1)/2}}=\eta_1$ 
and so $E(\eta_1)^{1-\mu_1^{(p-1)/2}}=E(1)$, 
or $E(\eta_1)^2=E(1)$, so $E(\eta_1)^2$ is of rank null and  therefore $E(\eta_1)= E(\eta_1^2)^{(p+1)/2}$ is also of rank null, contradiction. The same for $\mu_i,\quad i=1,\dots, r_\beta$.
\end{proof}
\end{lem}
%
\subsection{$\pi$-adic congruences on $p$-unit group $F$ of $\Q(\zeta)$}
The results on structure of relative $p$-class group $C_p^-$  of subsection
\ref{s108311} p. \pageref{s108311} can be translated to some results on structure of the group $F$:
from $\eta_i^{p-1}\equiv 1\modu \pi$ and from  $<E(\eta_i^{p-1})>=<E(\eta_i)>$,  we can always, without loss of generality, choose the determination  $\eta_i$ such that $\eta_i\equiv 1\modu\pi$.
We have proved that 
\begin{equation}\label{e203082}
\begin{split}
& \eta_i\equiv 1\modu\pi,\\
& \sigma(\eta_i)\equiv \eta_i^{\mu_i}\modu\pi^{p+1}.
\end{split}
\end{equation}
Then, starting of this relation (\ref{e203082}), similarly to lemma  \ref{l108171} p. \pageref{l108171} we get:
%
\paragraph{$\pi$-adic congruences of unit group $F=\{\Z[\zeta+\zeta^{-1}]^*/(\Z[\zeta+\zeta^{-1}]^*)^p\}/<-1>$}
This theorem summarize our $\pi$-adic approach on group of $p$-units $F$.
\begin{thm} { *** }\label{t207271}
With a certain ordering of index $i=1,\dots,\frac{p-3}{2}$, 
there exists a fundamental system of units $\eta_i,\quad i=1,\dots,\frac{p-3}{2}$, of the group
$F=\{\Z[\zeta+\zeta^{-1}]^*/(\Z[\zeta+\zeta^{-1}]^*)^p\}/<-1>$ verifying the relations:
\begin{equation}\label{e201274}
\begin{split}
& \eta_i\in\Z[\zeta+\zeta^{-1}]^*,\quad i=1,\dots,\frac{p-3}{2},\\
& \mu_i=u_{2n_i},\quad 1\leq n_i\leq\frac{p-3}{2},\quad i=1,\dots,\frac{p-3}{2},\\
& \sigma(\eta_i)=\eta_i^{\mu_i}\times\varepsilon_i^p,
\quad \varepsilon_i\in\Z[\zeta+\zeta^{-1}]^*,
\quad i=1,\dots,\frac{p-3}{2},\\
& \sigma(\eta_i)\equiv \eta_i^{\mu_i}\modu\pi^{p+1},\quad i=1,\dots,\frac{p-3}{2},\\
& \pi^{2n_i} \ \| \ \eta_i-1,\quad i=1,\dots,r_p^+,\quad \eta_i \mbox{\ not primary},\\
& \pi^{a_i(p-1)+2n_i} \ \| \ \eta_i-1,\quad a_i\in\N,\quad a_i>0,
\quad i=r_p^++1,\dots,r_p^-,\quad \eta_i \mbox{\  primary}.\\
& \pi^{a_i(p-1)+2n_i} \ \| \ \eta_i-1,\quad a_i\in\N,\quad a_i\geq 0,
\quad i=r_p^-+1,\dots,r_p,\quad \eta_i \mbox{ primary or not primary} \\
& \pi^{2n_i} \ \| \ \eta_i-1,\quad i=r_p+1,\dots,\frac{p-3}{2},\quad \eta_i \mbox{\ not primary.}\\
\end{split}
\end{equation}
\begin{proof}$ $
\begin{enumerate}
\item
We are applying  in this  situation  the same $\pi$-adic theory to $p$-group of units $F=\Z[\zeta+\zeta^{-1}]^*/(\Z[\zeta+\zeta^{-1}]^*)^p$ than to relative $p$-class group $C_p^-$ in subsection \ref{s108311} p. \pageref{s108311}, with a supplementary result for units due to Denes, see
Denes \cite{de1} and \cite{de2} and Ribenboim \cite{rib} (8D) p. 192.
\item
Similarly to decomposition of components of $C_p$ in singular primary and singular not primary components, 
the rank $\frac{p-3}{2}$ of $F$  has two components $\rho_1$ and $\frac{p-3}{2}-\rho_1$ where $\rho_1$ corresponds to the maximal number of independant units $\eta_i$  primary and 
$\rho_2=\frac{p-3}{2}-\rho_1$ to the units $\eta_i$  not primary.
\end{enumerate}
\end{proof}
\end{thm}
%

%
%
The next lemma for the unit group $\Z[\zeta+\zeta^{-1}]^*$ is the translation of similar lemma \ref{l203031} p. \pageref{l203031} for the relative $p$-class group $C_p^-$.
\begin{lem}{ *** }\label{l203042}
Let $\eta_{1}, \eta_{2}$ defined by relation (\ref{e201274}) p. \pageref{e201274}.
If $\mu_{1}=\mu_{2}$ then $\eta_{1}\times\eta_{2}^{-1}$ is  a primary unit.
\end{lem}
%
The group $F=\Z[\zeta+\zeta^{-1}]^*/(\Z[\zeta+\zeta^{-1}]^*)^p$ can be written as  the direct sum 
$F=F_1\oplus F_2$  of a subgroup $F_1$ with $\rho_1$  primary units ($p$-rank $\rho_1$ of $F_1$) and of a subgroup $F_2$ with $\rho_2=\frac{p-3}{2}-\rho_1$ fundamental not primary units 
($p$-rank $\rho_2$ of $F_2$): towards this assertion, observe that if $\eta_1$ and $\eta_2$ are two 
not primary units with $\sigma(\eta_1)\times \eta_1^{-\mu}\in(\Z[\zeta+\zeta^{-1}]^*)^p$ and 
$\sigma(\eta_2)\times \eta_2^{-\mu}\in(\Z[\zeta+\zeta^{-1}]^*)^p$  then $\eta_1\times\eta_2^{-1}$ is a primary unit and it always possible to replace $\{\eta_1,\eta_2\}$ by 
$\{\eta_1\times\eta_2^{-1},\eta_2\}$ in the basis of $F$, so to {\it push} all the primary units in $F_1$ and to make the set $F_2$ of not primary units as a group.
Observe that $\rho_1$ can be seen  also as  the maximal number of independant  primary units in $F$.
%
\paragraph{Structure theorem of unit group $F=\Z[\zeta+\zeta^{-1}]^*/(\Z[\zeta+\zeta^{-1}]^* )^p$}
\begin{thm}{ *** }\label{t203041}

Let $r_p^-$ be the relative $p$-class group of $\Q(\zeta)$.
Let $r_p^+$ be the $p$-class group of $\Q(\zeta+\zeta^{-1})$.
Let $\rho_1$ be the number of independant primary units of $F$. Then
\begin{equation}\label{e203072}
 r_p^--r_p^+\leq \rho_1\leq   r_p^-.
\end{equation}
\begin{proof}$ $
We apply Hilbert class field theory: for a certain order of the indexing of $i=1,\dots,\frac{p-3}{2}$:
\begin{enumerate}
\item
There are exactly $r_p^+$ independant unramified cyclic extensions
\begin{displaymath} 
\Q(\zeta,\omega_i)/\Q(\zeta), \quad \omega_i^p\in\Z[\zeta+\zeta^{-1}]-\Z[\zeta+\zeta^{-1}]^*
\quad i=1,\dots,r_p^+.
\end{displaymath} 
\item
There are exactly
$r_p^--r_p^+=r_p^-$ independant unramified cyclic extensions
\begin{displaymath} 
\Q(\zeta,\omega_i)/\Q(\zeta),\quad  \omega_i^p\in\Z[\zeta+\zeta^{-1}]^*,
\quad i=r_p^++1,\dots,r_p^-.
\end{displaymath}
\item
There is a number $n$ on independant unramified cyclic extensions witn $0\leq n\leq r_p^+$ with 
\begin{displaymath} 
\Q(\zeta,\omega_i)/\Q(\zeta),\quad  \omega_i^p\in\Z[\zeta+\zeta^{-1}]^*,
\quad i=r_p^-+1,\dots,r_p.
\end{displaymath}
\item
There are no  independant unramified cyclic extensions  with 
\begin{displaymath} 
\Q(\zeta,\omega_i)/\Q(\zeta),\quad  \omega_i^p\in\Z[\zeta+\zeta^{-1}]^*,
\quad i=r_p+1,\dots,\frac{p-3}{2}.
\end{displaymath}
\end{enumerate}
\end{proof}
\end{thm}
%
\paragraph{The case $\mu=u_{2n}$ with  $ 2n>\frac{p-1}{2}$}\label{s201111}
\paragraph{ }
In the  next theorem  we shall investigate more deeply the consequences of the congruence $\eta_i\equiv 1\modu \pi^{2n_i}$  when $2n_i> \frac{p-1}{2}$. We give an explicit congruence formula in that case.
To simplify notations, we take $\eta,\mu,n$ for $\eta_i,\mu_i, n_i$.
%
%
The next theorem for the $p$-unit group $F=\{\Z[\zeta+\zeta^{-1}]^*/(\Z[\zeta+\zeta^{-1}]^*)^p\}/<-1>$ is the translation of similar theorem \ref{l202211} p. \pageref{l202211} for the relative $p$-class group $C_p^-$.
\begin{thm}{ *** }\label{l203044}
Let $\mu=u_{2n},\quad p-3\geq 2n>\frac{p-1}{2}$, corresponding to 
$\eta\in\Z[\zeta+\zeta^{-1}]^*$ defined in relation (\ref{e201274}) p. \pageref{e201274}, so 
$\sigma(\eta)= \eta^{\mu}\times\varepsilon^p,\quad\varepsilon\in\Z[\zeta+\zeta^{-1}]^*$.
Then $\eta$ verifies the explicit formula: 
\begin{equation}\label{e203043}
\eta\equiv 1-\frac{\gamma_{p-3}}{\mu-1}\times
(\zeta+\mu^{-1}\zeta^u+\dots+\mu^{-(p-2)}\zeta^{u_{p-2}}) \modu \pi^{p-1},\quad \gamma_{p-3}\in\Z.
\end{equation}
\end{thm}
%
%
%
%
%

Roland Qu\^eme

13 avenue du ch\^ateau d'eau

31490 Brax

France

2005 fev 04

tel : 0033561067020

mailto: roland.queme@wanadoo.fr

home page: http://roland.queme.free.fr/

\end{document}